\documentstyle{amsppt}
\loadbold
\parskip=5.0pt
\topmatter
\title  Global existence of solutions of the
Liquid Crystal  flow for  the Oseen-Frank model
\endtitle
\rightheadtext{The Liquid Crystal  flow} \leftheadtext{M.-C. Hong
and Z. Xin}
\author {Min-Chun Hong and Zhouping Xin}\endauthor
  \address {Min-Chun Hong, Department of Mathematics,
  The University of Queensland, Brisbane, QLD 4072, Australia.
 Email:  hong\@maths.uq.edu.au
  \medskip
  Zhouping Xin, The Institute of Mathematical Sciences,
The Chinese University of Hong Kong,  Shatin, N.T., Hong Kong.
 }
\endaddress

\abstract {In the first part of this paper,  we establish global
existence of solutions of the liquid crystal (gradient) flow for
the well-known Oseen-Frank model. The liquid crystal  flow  is a
prototype of equations from the Ericksen-Leslie system in the
hydrodynamic theory  and generalizes the heat flow for harmonic
maps  into the $2$-sphere. The Ericksen-Leslie system is a system
of the Navier-Stokes equations coupled with the liquid crystal
flow. In the second part of this paper, we also prove global
existence of solutions of the Ericksen-Leslie system for a general
Oseen-Frank model in $\Bbb R^2$.}
\endabstract

\subjclass{AMS 35K50,  35Q30}\endsubjclass

\keywords {Liquid crystals flow, Navier-Stokes equations
}\endkeywords
 \endtopmatter
\document

\def\R{{\Bbb R}} 

  \def\laplacian{\bigtriangleup}
  \def\a{\alpha}
\def\intave#1{-\kern-10.7pt\int_{\,#1}}
\def\b{\beta}
\def\D{\nabla}
\def\g{\gamma}

\def\<{\langle}
\def\>{\rangle}
\def\({\left(}
\def\){\right)}
\def\esssup{\operatornamewithlimits{ess\,sup}}
\def\limsup{\operatornamewithlimits{lim\,sup}}
\def\intave#1{-\kern-10.7pt\int_{\,#1}}

\head {\bf 1. Introduction}\endhead

A {\it liquid crystal} is a state of matter  intermediate between
a crystalline solid and a normal isotropic liquid. Research into
liquid crystals is an area of
  a very successful synergy between mathematics and
physics.  There are a lot of analytical and computational issues,
which arise in the attempt to study static equilibrium
configurations. Numerical and experimental analysis has shown that
equilibrium configurations are expected to have point and line
singularities (\cite {K}). Mathematically, Hardt, Kinderlehrer and
Lin  in their fundamental papers \cite {HKL1} and \cite {HKL2}
proved the existence of an energy minimizer $u$ of the liquid
crystal functional and showed that a minimizer $u$ is smooth away
from a closed set $\Sigma$ of $\Omega$. Moreover, $\Sigma$ has
Hausdorff dimension strictly less than one. In \cite {AL}, Almgren
and Lieb did some related analysis indicating that the phenomenon
is of wider interest. In physical theory, an equilibrium
configuration corresponds to a critical point, not necessarily an
energy minimizer, of the liquid crystal energy. Critical points
are much harder to understand mathematically than minima.   From
the above result of Hardt, Kinderlehrer and Lin, minimizers cannot
have line singularities. Following the work of
Bethuel-Brezis-Coron on harmonic maps in \cite {BBC}, Giaquinta,
Modica and Soucek \cite {GMS2} found a relaxed energy for the
liquid crystal systems, whose minimizers are also equilibrium
configurations. On the other hand, Giaquinta, Modica and Soucek
\cite {GMS1} also proved that  minimizers of the relaxed energy
for harmonic maps are smooth away from a $1$-dimensional singular
set. Further developments on the regularity results on harmonic
maps were surveyed in \cite {GMS3}. There is an interesting open
problem to prove that minimizers of the relaxed liquid crystal
energy have line singularities. The first author in \cite {Ho3}
proved partial regularity of minimizers  of the modified relaxed
energy of the liquid crystal energy. However, the partial
regularity of minimizers of the relaxed energy for liquid crystals
is still mysterious. In some related studies of liquid crystals,
Bauman, Calderer, Liu and Phillis  \cite {BCCP} studied the
Landau-de Gennes free energy used to describe the transition
between chiral nematic and smectic liquid crystal phase, Lin and
Pan \cite {LP} used the Landau-de Gennes models to investigate the
magnetic field induced instabilities in liquid crystals, and the
existence of infinite many liquid crystal equilibrium
configurations prescribing the same boundary was obtained in \cite
{Ho2}.

A general description of the static theory  of liquid crystals is
given by Ericksen in \cite {Er}. A liquid crystal is composed of
rod like molecules which display orientational order, unlike a
liquid, but lacking the lattice structure of a solid.   The
kinematic variable in the nematic and cholesteric phase may be
taken to the optic axis, which is a unit vector field $u$ in a
region $\Omega\subset\R^3$ occupied by the materials. The liquid
crystal energy for a configuration $u\in H^1(\Omega ; S^2)$ is
given by
$$E(u; \Omega )=\int_{\Omega} W(u, \nabla u)\,dx,\tag 1.1$$
where the Oseen-Frank density $W (u,\nabla u)$, depending  on
positive material constants $k_1$, $k_2$, $k_3$ and $k_4$, is
given by
$$W(u,\nabla u)=k_1(\text{div } u)^2+k_2 (u\cdot\text{curl }u)^2+k_3
|u\times \text{curl } u|^2+k_4[\text {tr} (\nabla u)^2-(\text {div
}u)^2].
$$
 Without loss of generality, as in \cite {HKL1} or \cite
{GMS3}, we rewrite the density
$$W(u,\nabla u) = a|\nabla u|^2 +
V(u, \nabla u),\quad a =\min\{k_1, k_2,k_3\}>0,\tag 1.2$$ where
$$V(u, \nabla u)=(k_1-a)(\text{div }
u)^2+(k_2-a) (u\cdot\text{curl }u)^2+(k_3-a) |u\times \text{curl }
u|^2. $$

A static equilibrium configuration corresponds to an extremal
(critical point) of the energy functional $E$ in $H^1(\Omega
,S^2)$. The Euler-Lagrange system for the general Oseen-Frank
functional (1.1) (see details in Appendix) is:
$$\split
&\quad \nabla_{\a} \left [W_{p_{\a}^i} (u,\nabla u)- u^lu^i
V_{p_{\a}^l}(u,\nabla u) \right ]-W_{u^i}(u,\nabla
u)+W_{u^l}(u,\nabla u)u^l u^i\\
&+W_{p_{\a}^l}(u,\nabla u)\nabla_{\a}u^l u^i
+V_{p_{\a}^l}(u,\nabla u) u^l\nabla_{\a}u^i=0\quad \text {in
}\Omega
\endsplit\tag 1.3
$$
for $i=1,2,3$, where we adopt the standard summation convention.
In a special case of $k_1=k_2=k_3$, the system (1.3) becomes the
harmonic map equations into $S^2$. However, the equilibrium system
associated to the energy functional (1.1)
 is not elliptic for every choice of the constants $k_1$, $k_2$ and $k_3$.

In the first part of this paper, we investigate the  liquid
crystal flow for a model with the Oseen-Frank  density (1.2). For
a domain $\Omega$ in $\R^3$ or in $\R^2$, a map
$u(x,t):\Omega\times [0,\infty )\to S^2$ is a solution of the
liquid crystal flow if $u$ satisfies
$$\split
\frac {\partial u^i}{\partial t}=&\nabla_{\a} \left [W_{p_{\a}^i}
(u,\nabla u)- u^lu^i V_{p_{\a}^l}(u,\nabla u) \right
]-W_{u^i}(u,\nabla
u)\\
&+W_{u^l}(u,\nabla u)u^l u^i+W_{p_{\a}^l}(u,\nabla
u)\nabla_{\a}u^l u^i +V_{p_{\a}^l}(u,\nabla u) u^l\nabla_{\a}u^i
\endsplit \tag 1.4
$$ in $\Omega\times  [0,\infty )$ for $i=1,2,3$.

 The  flow  equation (1.4) is a prototype of equations from the
Ericksen-Leslie system in the hydrodynamic theory (cf. \cite
{Er}). The liquid crystal flow (1.4) also generalizes the heat
flow for harmonic maps  into the $2$-sphere. Since the seminal
work of Eells-Sampson \cite {ES}, many  work on the heat flow for
harmonic maps have been made. In 2 dimensional case, Struwe
\cite{St1} established  global existence of the weak solution of
the harmonic maps flow with initial data, where the solution is
smooth except for a finite number of singularities. In higher
dimensional cases, Chen-Struwe \cite{CS} proved global existence
of partially regular solutions to the harmonic map flow. Since
(1.4) is not parabolic, the system of the liquid crystal flow is
complicated, so the question on global existence  for the liquid
crystal flow (1.4) for the Oseen-Frank model remains unresolved.
In this paper, we prove global existence of solutions of the
liquid crystal flow in 2D.

We set
$$H_{b}^{1} (\R^2; S^2):=\left \{u :  \quad u-b\in H^1(\R^2; \R^3), \quad |u|=1\, \text { a.e. in }
\R^2 \right \}
$$
for a constant vector $b\in S^2$.

Then, one of our main results in this paper is the following
global existence for this flow in 2D (i.e. $u$ is a constant along
a direction in $\R^3$):

 \proclaim
{Theorem A} Let $u_0\in H_b^{1} (\R^2; S^2)$  be a given map. Then
there exists a global weak solution $u(x,t): \R^2\times
[0,+\infty)\to S^2$ of (1.4) with initial value $u(0)=u_0$ such
that $u$ is smooth in $\Omega\times [0,+\infty )$ except for a
finite number of singularities $\{(x^l_i,T_l)\}_{l=1}^K\in
\R^2\times [0,+\infty)$ with an integer $K>0$ depending on $u_0$.
Moreover, there are two constants $\varepsilon_0>0$ and $R_0>0$
such that each singular point $\ x_{i}^{l} $ at the time $T_l$ is
characterized by the condition
$$\limsup_{t \nearrow  T_{l}} E\left( u\left(
t\right)  ,B_{R}\left( x_{i}^{l}\right) \right) \geq
\varepsilon_{0}$$
 for any  $R>0$ with $R\leq R_0$.
\endproclaim
This result can be regarded as an extension of the well-known
result of Struwe in \cite {St1} on the heat flow for harmonic maps
in dimension two. Since the liquid crystal flow is not a parabolic
system, the flow (1.4) is more complicated than the harmonic map
flow. In particular, we can not apply the well-known theory  of
partial differential equations  directly to prove the local
existence for the liquid crystal flow. Instead, we consider a
family of Ginzburg-Landau approximation flows to prove local
existence of solutions to (1.4). To prove Theorem A, we  need to
get a $L^2$-estimate of $\nabla^2u$ similarly to  one in \cite
{St1}. However, the  flow (1.4) is not a parabolic system, so we
overcome the difficulties due to the term $\nabla_{\a}  [ u^lu^i
V_{p_{\a}^l}(u,\nabla u)]$ by using the fact that $|u|=1$ as
observed in  \cite {Ho1}.

In the second part of this paper, we investigate the
Ericksen-Leslie system with the   Oseen-Frank density $W(u,\nabla
u)$ in (1.2). In the 1960's, Ericksen \cite {Er} and Leslie \cite
{Le} established the hydrodynamic theory of liquid crystals
independently. The Ericksen-Leslie theory describes the dynamic
flow of liquid crystals, including the velocity vector $v$ and
direction vector $u$ of the fluid. Let $v=(v^1,v^2,v^3)$ be the
velocity vector of the fluid and $u=(u^1,u^2,u^3)$ the unit
direction vector. The Ericksen-Leslie system in $\Omega\times
[0,\infty )$ is given by (e.g. \cite {L1} and \cite {LL1})
$$v^i_t+(v\cdot\nabla )v^i-\nu\laplacian v^i+\nabla_{x_i} P=-\lambda
\nabla_{x_j}(\nabla_{x_i}u^kW_{p_j^k}(u, \nabla u)),\tag 1.5
$$
$$\nabla \cdot v=0,\tag 1.6
$$
$$\split u^i_t+(v\cdot \nabla  )u^i =&  \nabla_{\a} \left [W_{p_{\a}^i}
(u,\nabla u)- u^ku^i V_{p_{\a}^k}(u,\nabla u) \right
]-W_{u^i}(u,\nabla
u)\\
&+W_{u^k}(u,\nabla u)u^k u^i+W_{p_{\a}^l}(u,\nabla
u)\nabla_{\a}u^l u^i +V_{p_{\a}^k}(u,\nabla u) u^k\nabla_{\a}u^i
\endsplit \tag 1.7
$$
for $i=1,2,3$,  prescribing the boundary condition
$$v(x,t)=0, \quad u(x,t)=u_0(x), \quad \forall (x,t)\in\partial \Omega
\times (0,\infty ) \tag 1.8
$$
and with initial data
$$v(x,0)=v_0(x), \quad u(x,0)=u_0(x),  \quad \text {div } v_0=0\quad \forall x\in \Omega . \tag
1.9
$$
Here $\nu$, $\lambda$ are given positive
constants, and $P$ is the pressure.

The system (1.5)-(1.7) is a system of the  Navier-Stokes equations
coupled with the liquid crystal flow (1.4). The study of the
 Navier-Stokes equations is of great interest.  Tremendous results on the
existence and partial regularity  for the Navier-Stokes equations
have been established (e.g. \cite {Sc}, \cite {CKN}, \cite {L2},
\cite {TX}). In this paper, we are only concentrating on the
existence of solutions of the  Ericksen-Leslie system.
 Since the  functional
$E(u; \Omega)$ in (1.1)   with the constraint $|u|=1$ is
complicated, one  considers Ginzburg-Landau functionals
 $$E_{\varepsilon}(u; \Omega )=\int_{\Omega} \left [ W(u, \nabla u)+ \frac 1 {2\varepsilon^2}
 (1-|u|^2)^2\,\right ]\,dx
 $$
 for any function $u\in H^1(\Omega ;\R^3)$. Then, the approximating Ericksen-Leslie system
 is given by
$$v^i_t+(v\cdot\nabla )v^i-\nu\laplacian v^i+\nabla_{x_i} P=
-\lambda  \nabla_{x_j}(\nabla_{x_i}u^kW_{p_j^k}(u,\nabla u)),
\quad  \tag 1.10
$$
$$ \nabla \cdot v=0,\tag 1.11
$$
$$
u^i_t+(v\cdot \nabla  )u^i =   \nabla_{\a} \left [W_{p_{\a}^i} (u,
\nabla u) \right ]-W_{u^i}(u,\nabla u)+\frac 1 {\varepsilon^2} u^i
(1-|u|^2) \tag 1.12
$$
for  $i=1,2,3$,  prescribing the boundary condition (1.8) and
initial condition
 (1.9).

In the case of $k_1=k_2=k_3$, Lin and Liu  \cite {LL1} proved
global existence of the the classical solution of (1.10)-(1.12)
with (1.8)-(1.9) in dimension two and   the weak solution of the
same system  in dimension three. Lin and Liu in \cite {LL2} also
analyzed the limit of solutions
$(v_{\varepsilon},u_{\varepsilon})$ of (1.10)-(1.12) as
$\varepsilon\to 0$, but it is not clear that the limiting solution
satisfies the original Ericksen-Leslie system (1.5)-(1.7) with
$|u|=1$. Therefore, there is an interesting question to establish
the global existence of solutions of (1.5)-(1.7) with (1.8)-(1.9).
The question for the case of $k_1=k_2=k_3$ has been answered by
the first author in \cite {Ho3} in $\R^2$  and Lin-Lin-Wang \cite
{LLW} in a general case for a domain of $\R^2$ independently. The
system (1.5)-(1.7) or (1.10)-(1.12) for the general Oseen-Frank
model is
 more complicated than the system for the case of $k_1=k_2=k_3$ since there is no maximum principle for the parabolic
system (1.12) in the case $k_1\neq k_2$ (see \cite {A}) and  the
term $W_{u^i} (u,\nabla u)$ in (1.12) will cause a trouble to
prove  global existence for the system.

In this paper, we will prove global existence of weak solutions to
the Ericksen-Leslie system (1.5)-(1.7) for a general Oseen-Frank
model in $\R^2$. More precisely, we have

  \proclaim
{Theorem B}  Let $(u_0, v_0)\in H_b^{1} (\R^2; S^2)\times
L^2(\R^2,\R^2)$  be  given  initial data with $\text{dvi } v_0=0$.
Then, there exists a global weak solution $(u, v):\R^2\times
[0,+\infty)\to S^2\times \R^2$ of (1.5)-(1.7) with initial values
(1.9), where the solution $(u,v)$ is smooth in $\R^2\times
((0,+\infty )\backslash \{T_l\}_{l=1}^L)$ for a finite number of
times $\{T_l\}_{l=1}^L$. Moreover, there are two constants
$\varepsilon_0>0$ and $R_0>0$ such that each singular point
$\left( x_{i}^{l},T_{l}\right) \in\Sigma \times \{T_l\}$ is
characterized by the condition
$$\limsup_{t \nearrow  T_{l}} \int_{B_{R}( x_{i}^{l})} |\nabla
u\left(\cdot , t\right)|^2 +|v\left(\cdot , t\right)|^2  \,dx \geq
\varepsilon_{0}$$
 for any  $R>0$ with $R\leq R_0$.
\endproclaim

The main idea  to prove Theorem B is to combine the idea in \cite
{Ho3} with the proofs of Theorem A.  The first key step is to
prove local existence of solutions of the system (1.5)-(1.7) by
considering the approximation system (1.10)-(1.12). To prove
global existence of solutions to (1.5)-(1.7), one of key steps is
to get a $L^2$-estimate of $\nabla^2u$ and $\nabla v$ in
$\R^2\times [0,T]$ under  a small energy condition as in \cite
{St1}. To show the regularity of the weak solution $(u,v)$ of
(1.5)-(1.7) in $\R^2\times (0,T)$, we establish a local energy
inequality under the small energy condition, which was first used
by Struwe in \cite {St2} for the $H$-system flow. Finally, we
prove regularity of solutions by controlling $L^2$-estimate of
$\nabla^2u$ and $\nabla v$ in $\R^2$ for $t\in (0,T)$. Since (1.7)
is not a parabolic system, the proof of Theorem B is more
difficult than one for the case of $k_1=k_2=k_3$ in \cite {Ho3}.
We overcome a number of difficulties  on the regularity and
uniqueness for the systems by employing the invariance of the
density (1.2) after a rotation.

The rest of the paper is organized as follows. In Section  2, we
prove the global existence for the liquid crystal flow in 2D. Some
global estimates for (1.5)-(1.7) are established in Section 3.
Then, we complete a proof of Theorem B in Section 4. Finally, the
regularity issue for the systems is dealt in Section 5.

\head {\bf 2. Existence of partial regular solutions of the liquid
crystal flow}\endhead In this section, we consider the flow (1.4)
in $\R^2$. For simplicity of notations, $u$ is assumed to be a
constant along $x_3$-direction in $\R^3$; i.e. $\frac {\partial
u}{\partial x_3}=0$.

For any two positive constants $\tau$ and $T$ with $\tau < T$, we
define
$$\split
 V(\tau, T):= \{ u: \R^2\times &[\tau , T]\to S^2\, |\quad  u \quad \text
  {is measureable and satisfies}\\
&
 \esssup_{\tau\leq t\leq T} \int_{\R^2} |\nabla u(\cdot ,t)|^2\, dx
 +\int_{\tau}^T \int_{\R^2} |\nabla^2 u|^2 +|\partial_t u|^2 \,dx\,dt <\infty
  \}.
 \endsplit
 $$

\proclaim {Lemma 1} Let $u\in V(0 ,T)$ be a solution  of the
system (1.4) with initial value $u_0\in H^1(\R^2,S^2)$. Then, for
any $t_1\in [0,T]$
$$\int_{\R^2\times (0 ,t_1)}|\partial_t u|^2\,dx\,dt + E(u(t_1))\leq
E(u_0).\tag 2.1
$$
Moreover,  for all $t\in [0,T]$,  $x_0\in\R^2$ and $R>0$, it holds
that
$$ \split
&\quad \int_{B_R(x_0)}W(u (x,t),\nabla u(x,t))\,dx\\
&\leq \int_{B_{2R}(x_0)}W(u_0(x),\nabla u_0(x))\,dx +C\frac
t{R^2}\int_{\R^2} |\nabla u_0|^2\,dx, \endsplit\tag 2.2
$$
where $C$ is a constant.
\endproclaim
\demo{Proof} Multiplying (1.4) by $\frac {\partial u^i}{\partial
t}$ yields
$$\int_{\R^2}|\frac {\partial u}{\partial t} |^2\,dx
=-\int_{\R^2} W_{p_{\a}^i} (u,\nabla u) \frac {d}{dt} \nabla_{\a}
u^i\,dx- \int_{\R^2}W_{u^i}(u,\nabla u) \frac {\partial
u^i}{\partial t}\,dx.
$$
This implies
$$\int_{\R^2}|\frac {\partial u}{\partial t} |^2\,dx
+\frac {d}{dt} \int_{\R^2} W(u,\nabla u)\,dx=0.
$$
(2.1) follows from integrating the above identity.

Let $\phi\in C^{\infty}_0(B_{2R}(x_0))$ be a cut-off function
satisfying $0\leq \phi\leq 1$, $|\nabla \phi |\leq C/R$ and
$\phi\equiv 1$ on $B_R(x_0)$. Multiplying (1.4) by $\frac
{\partial u^i}{\partial t}\phi^2$ and then using Young's
inequality yields
$$\split
&\int_{\R^2}|\frac {\partial u}{\partial t} |^2\phi^2\,dx +\frac d
{dt}\int_{\R^2} W (u(x,t),\nabla u(x,t))\phi^2 \,dx\leq
C\int_{\R^2}|\frac {\partial u}{\partial t} | \,|\nabla u|
|\phi\nabla
\phi |\,dx\\
&\leq \frac 1 2\int_{\R^2}|\frac {\partial u}{\partial t}
|^2\phi^2\,dx +C\int_{\R^2} W(u,\nabla u)|\nabla \phi|^2\,dx.
\endsplit
$$
Then,  (2.2) follows from using (2.1) and integrating the above
inequality.\qed
\enddemo
It follows from \cite {St1} that
 \proclaim {Lemma 2} There are  constants $C$ and $R_0$
such that for any $u\in V(0 ,T)$ and any $R\in (0, R_0]$,  we have
$$\split  \quad \int_{\R^2\times [0 ,T]} |\nabla u|^4\,dx\,dt\leq &C\esssup_{0\leq t\leq T, x\in \R^2} \int_{B_R(x)}|\nabla
u(\cdot ,t)|^2\,dx  \\& \cdot  (\int_{\R^2\times [0 ,T]} |\nabla^2
u|^2  \,dx\,dt +R^{-2} \int_{\R^2\times [0 ,T]} |\nabla u|^2
\,dx\,dt ).
\endsplit
$$
\endproclaim

\proclaim {Lemma 3} Let $u\in V(0 ,T)$ be a   solution   of (1.4)
with initial smooth value $u_0\in H^1$. Then there are constants
$\varepsilon_1$ and $R_0>0$ such that if
$$\esssup_{0\leq t\leq T, x\in \R^2} \int_{B_R(x)}|\nabla
u(\cdot ,t)|^2\,dx<\varepsilon_1
 $$ for  any
$R\in (0, R_0]$, then
$$\int_{\R^2\times [0 ,T]} |\nabla^2 u|^2 \,dx\,dt\leq C
E(u_0) \, (1 +T R^{-2}),\tag 2.3
$$
$$\int_{\R^2\times [0 ,T]} |\nabla u|^4 \,dx\,dt\leq C
\varepsilon_1 E(u_0) \, (1 +T R^{-2}).\tag 2.4
$$
\endproclaim
\demo{Proof} Multiplying   (1.4) by $\laplacian u^i$ yields
$$\split
&\quad \int_{\R^2}\frac {\partial u^i}{\partial t}\laplacian
u^i\,dx\\
&=\int_{\R^2}  \nabla_{\a}  \left [W_{p_{\a}^i} (u,\nabla
u)- u^ku^i V_{p_{\a}^k}(u,\nabla u) \right ]\laplacian u^i\,dx\\
& -\int_{\R^2}W_{u^i}(u,\nabla u)(\laplacian u^i -u^k
u^i\laplacian u^k)\,dx +\int_{\R^2}W_{p_{\a}^k}(u,\nabla
u)\nabla_{\a}u^k u^i \laplacian
u^i \,dx\\
&+\int_{\R^2}V_{p_{\a}^k}(u,\nabla u) u^k\nabla_{\a}u^i \laplacian
u^i\,dx:=I_1+I_2+I_3+I_4.
\endsplit
$$
Note that the terms $I_2$ and $I_3$ of the above identity can be
controlled by $C|\nabla u|^2|\laplacian u|$. It suffices to
estimate terms $I_1$ and $I_4$. Since $|u|^2=1$, $-u^i\laplacian
u^i =|\nabla u|^2$. We note
$$\split \nabla_{\a}[u^ku^i V_{p_{\a}^k}(u,\nabla u) ]=&\nabla_{\a} u^k
u^iV_{p_{\a}^k}(u,\nabla u)+  u^k \nabla_{\a}
u^iV_{p_{\a}^k}(u,\nabla u)\\
&+ u^k u^i \nabla_{\a} V_{p_{\a}^k}(u,\nabla u).
\endsplit
$$
Integration by parts twice yields
$$\split
I_1+I_4&=\int_{\R^2} \nabla_{\b}  \left  [W_{p_{\a}^i} (u,\nabla
u) \right ]\nabla^2_{\a\b}u^i\,dx+ \int_{\R^2} \nabla_{\a} u^k
u^iV_{p_{\a}^k}(u,\nabla u) \laplacian u^i\, dx\\
 &- \int_{\R^2}
u^k   \nabla_{\a} V_{p_{\a}^k}(u,\nabla u)|\nabla u|^2\,dx.
\endsplit
$$
Note
$$ \nabla_{\a} V_{p_{\a}^k}(u,\nabla u)=V_{p_{\a}^k p_{\g}^l}
(u,\nabla u) \nabla^2_{\g\a}u^l+V_{p_{\a}^ku^l}(u,\nabla
u)\nabla_{\a} u^l
$$
and
$$\nabla_{\b} W_{p_{\a}^i}(u,\nabla u)=W_{p_{\a}^i p_{\g}^j}
(u,\nabla u) \nabla^2_{\g\b}u^j+W_{p_{\a}^iu^j}(u,\nabla u)
\nabla_{\b} u^j.
$$ This implies
$$\split
&\frac {d}{dt}\int_{\R^2} |\nabla u|^2\,dx +\int_{\R^2}
W_{p_{\a}^i p_{\g}^j} (u,\nabla u)\nabla^2_{\a\b}u^i
\nabla^2_{\g\b}u^j \,dx\\
&\leq C\int_{\R^2} |\nabla u|^2 (|\nabla u|^2+|\nabla^2u|)\,dx.
\endsplit\tag 2.5
$$
Since $W(u,p)$ is convex in $p$, it satisfies the ellipticity
$$W_{p_{\a}^i p_{\g}^j}
(u,\nabla u)\nabla^2_{\a\b}u^i \nabla^2_{\g\b}u^j\geq a|\nabla^2
u|^2
$$
for the constant $a>0$.
 Then, choosing $\varepsilon_1>0$ to be sufficiently
 small and applying Lemma 2 lead to (2.3) and (2.4).
 \qed
\enddemo

\proclaim {Lemma 4} Let $u\in V(0,T)$ be a solution of (1.4) with
initial value $u_0\in H^1$. Assume that
$$\esssup_{0\leq t\leq T, x\in \R^2} \int_{B_R(x)}|\nabla
u(\cdot ,t)|^2\,dx<\varepsilon_1
$$ for any
$R\in (0, R_0]$. Let $\tau\in(0,T]$ be any constant. Then it holds
for all $t\in [\tau, T]$,
$$\int_{\R^2} |\nabla^2 u(x,t) |^2 \,dx\leq C_0,\tag 2.6
$$
with a uniform constant depending only on $\tau$, $T$, $R_0$, and
$E(u_0)$.
\endproclaim
\demo {Proof} The proof is similar to \cite {St1; Lemma 3.10}.
Using a proper cut-off function if necessary, we assume in the
following proof that $\int|\partial_t\nabla u|^2 (\cdot,t) ,dx$ is
finite.

Differentiate (1.4) with respect to $t$, multiply the resulting
identity by $\partial_t u^i$, and then integrate to obtain
$$\split
& \ \frac 12 \frac{d}{dt} \int_{\R^2} |\partial_t u|^2\,dx +
\int_{\R^2}W_{p_{\a}^i p_{\b}^j}(u,\nabla u) \nabla_{\a}\,\partial_t\,u^i\,\nabla_{\b}\,\partial_t\,u^j\,dx\\
\leq & \ C \int_{\R^2} (|\partial_t\,u|^2\,|\nabla u|^2+|\nabla
u|\,
|\partial_t\,u|\,|\nabla\partial_t\,u|]\,dx \\
&+\int_{\R^2}
\partial_t\,u^i\,\nabla_{\a} \left [ u^i\,u^k\,V_{p^k_{\a}} (u,\nabla
\partial_t\,u)\right ]\,dx.
\endsplit
$$
Due to the convexity of $W(u,p)$ in $p$, there exists a positive
constant $a>0$ such that
$$\int_{\R^2}W_{p_{\a}^i p_{\b}^j}(u,\nabla u) \nabla_{\a}\,\partial_t\,u^i\,\nabla_{\b}\,\partial_t\,u^j\,dx
\leq a \int_{\R^2} |\nabla\partial_t\,u(x,t)|^2\,dx.$$ Since
$|u|=1$, so $\sum_i\,\partial_t\,u^i\,u^i=0$. And hence,
$$\int_{\R^2} \partial_t
u^i \nabla_{\a} \left ( u^i\,u^k\,V_{p^k_{\a}} (u,\nabla
\partial_t\,u)\right
)\,dx\leq C \int_{\R^2} |\partial_t\,u|\,|\nabla
u|\,|\nabla\partial_t\,u|\,dx.
$$
It follows from these and Cauchy's inequality that
$$
\frac{1}{2} \frac{d}{dt} \int_{\R^2} |\partial_t\,u(x,t)|^2\,dx
+\frac{a}{2}\int_{\R^2} |\nabla \partial_t\,u(x,t)|^2\,dx \leq C
\int_{\R^2} |\partial_t\,u|^2\,|\nabla u|^2\,dx.\tag 2.7
$$
Note that
$$\split
& \quad \ C\int_{\R^2} |\partial_t\,u|^2\,|\nabla u|^2\,dx\\
& \leq C
\left( \int_{\R^2}
|\partial_t\,u|^4\,dx\right)^{\frac{1}{2}}\,\left(\int_{\R^2}
|\nabla u|^4\,dx\right)^{\frac{1}{2}}\\
& \leq C \left(\int_{\R^2}
|\partial_t\,u|^2\,(x,t)\,dx\right)^{\frac{1}{2}}\,\left(\int_{\R^2}
|\partial_t\,\nabla
u|^2\,(x,t)\right)^{\frac{1}{2}}\,\left(\int_{\R^2}|\nabla
u|^4\,dx\right)^{\frac{1}{2}}\\
& \leq \frac{a}{4}\int_{\R^2}|\nabla\partial_t\,u(x,t)|^2\,dx +
\left(C\int_{\R^2}|\nabla u(x,t)|^4\,dx\right)
\,\int_{\R^2}|\partial_t\,u(x,t)|^2\,dx.
\endsplit
$$
This, together with (2.7), yields that for all $t\in(0,T]$,
$$\split
& \ \frac{d}{dt}\int_{\R^2} |\partial_t\,u(\cdot,t)|^2\,dx +
\frac{a}{2}\int_{\R^2}|\nabla\partial_t\,u(\cdot,t)|^2\,dx\\
\leq & \left(C\int_{\R^2}|\nabla
u(\cdot,t)|^4\,dx\right)\,\int_{\R^2}|\partial_t\,u(\cdot,t)|^2\,dt.
\endsplit\tag 2.8
$$
It follows from (2.8), (2.4), Lemma 1, and Gronwall's inequality
that for any $0<s\leq t\leq T$,
$$\split
\int_{\R^2}|\partial_t\,u(\cdot,t)|^2\,dt & \displaystyle \leq
\left( e^{C\int^t_s\!\!\int_{\R^2}|\nabla
u(\cdot,l)|^4\,dx\,dl}\right)\,\int_{\R^2}|\partial_t\,u(\cdot,s)|^2\,dx\\
& \leq
e^{C\varepsilon_1\,E(u_0)(1+TR^{-2})}\cdot\int_{\R^2}|\partial_t\,u(\cdot,s)|^2\,dx.
\endsplit
$$
Combining this with (2.1) shows that for any fixed $0<\tau<T$,
there exists a constant $C$ such that
$$\esssup_{\tau\leq t\leq T} \int_{\R^2} |\partial_t\,u (\cdot,t)|^2\,dt\leq
C\tau^{-1}\,E(u_0)\,e^{C\varepsilon_{1}\,E(u_0)(1+TR^{-2})},
$$
with a uniform constant $C$. On the other hand, using (2.5),
integration by parts yields that for any $t\in[\tau,T]$,
$$\split
& \ \int_{\R^2} |\nabla^2 u(\cdot,t)|^2\,dx \leq
C\int_{\R^2}|\nabla
u(\cdot,t)|^4\,dx +C\int_{\R^2}|\partial_t\,u(\cdot,t)|^2\,dx\\
\leq & \ C\varepsilon_1\int_{\R^2}|\nabla^2\,u(\cdot,t)|^2\,dx +
\frac{C\varepsilon_1}{R^2_0} E(u_0) +
C\int_{\R^2}|\partial_t\,u(\cdot,t)|^2\,dx.
\endsplit
$$
Combining this with (2.9) shows that for suitably small
$\varepsilon_1$, the desired estimate (2.6) holds with
$$C_0 \equiv CE(u_0)\left(\frac{\varepsilon_1}{R^2} + \tau^{-1}\,e^{C\varepsilon_1\,E(u_0)(1+TR^{-2})}
\right).\tag 2.10
$$

By the well-known Gagliardo-Nirenberg-Sobolev inequality, we have
for any $x\in\R^2$
$$\split
& |u(x,t_1)-u(x,t_2)|\leq C
\|u(x,t_1)-u(x,t_2)\|^{3/4}_{H^{2}( B_1(x))}\|u(x,t_1)-u(x,t_2)\|^{1/4}_{L^{2}(B_1(x))}\\
&\leq C (\sup_{\tau\leq t\leq T}\|\nabla
^2u(\cdot,t)\|^{3/4}_{L^{2}(\R^2)} +1)|t_1-t_2|^{1/8} \left
(\int_0^T \int_{\R^2}|\partial_t u|^2\,dx\,dt\right )^{1/8}\\
&\leq C|t_1-t_2|^{1/8}.
\endsplit
$$
It follows from (2.6) and Sobolev embedding theorem that $u(x,t)$
is
 H\"older continuous in $x$ uniformly for $t\in [\tau, T]$. Then we get
 that $u$ is H\"older continuous in  $C^{1/8}(\R^2\times [\tau, T])$ for any
 $T<T_1$. Due to Proposition 14 in Appendix, $u$ is in
 $C^{1,\frac 1 8}$. Hence, $u$ is regular in $(0,T_1)$.
\qed
\enddemo

\proclaim {Remark 5} Let $u\in V(0 ,T)$ be a   solution   of (1.4)
with initial value $u_0\in H_b^2$. Assume that there are constants
$\varepsilon_1$ and $R_0>0$ such that
$$\esssup_{0\leq t\leq T, x\in \R^2} \int_{B_R(x)}|\nabla
u(\cdot ,t)|^2\,dx<\varepsilon_1
 $$ for  any
$R\in (0, R_0]$.  Then, for any $t\in [0, T]$ and $R\leq R_0$, we
have
$$\int_{\R^2} |\nabla^2 u(x,t) |^2\,dx\leq C_1\equiv
C_1(||u_0||_{H^2_b}, C_0).
$$
\endproclaim

\proclaim {Theorem 6}  {\bf(Local existence)} For a map $u_0\in
H_b^1(\R^2, S^2)$, there is a solution $u\in V(0, t_1 )$ of (1.4)
with initial value $u_0$ for some $t_1>0$.
\endproclaim
\demo {Proof} For any map $u_0\in H_b^{1}(\R^2, S^2)$, it can be
approximated by a sequence of smooth maps in $H_b^{2}(\R^2, S^2)$.
Without loss of generality, we assume that $u_0\in H_b^{2}(\R^2,
S^2)$ is smooth. The liquid crystal flow is not a parabolic
system, so one can not apply the well-known local existence
theory. Instead, we prove the local existence by an approximation
of the Ginzburg-Landau flow in the following:
$$  \frac {\partial u_{\varepsilon}^i}{\partial
t}= \nabla_{\a} \left [W_{p_{\a}^i} (u_\varepsilon, \nabla
u_{\varepsilon}) \right ]-W_{u^i}(u_\varepsilon,\nabla
u_\varepsilon)+\frac 1 {\varepsilon^2} u_{\varepsilon}^i
(1-|u_{\varepsilon}|^2)
 \tag 2.11
$$
with  initial value $u_0\in H_b^2(\R^2, S^2)$ and $u_0\in
C^{\infty}$. Applying the standard local existence theory of
quasi-linear parabolic systems (cf. \cite {Ei} or \cite {Am}),
there is a local regular solution $u_{\varepsilon}$ of (2.11) with
initial value $u_{\varepsilon} (0)=u_0$.

For simplicity of notations, we define
$$\split
\tilde{V}(\tau,T)= & \ \biggl\{ u:\R^2\times[\tau,T]\rightarrow
\R^3| u
\text{\ is\ measurable\ and\ satisfies}\\
& \ \esssup_{\tau\leq t\leq T}\int_{\R^2}|\nabla u(\cdot,t)|^2\,dx
+ \int^T_\tau\!\!\int_{\R^2}(|\nabla^2
u|^2+|\partial_t\,u|^2)\,dx\,dt<\infty\biggl\},\\
& \ \ 0\leq\tau<T<+\infty,\\
e_{\varepsilon}(u)=&\ W(u,\nabla u)+\frac 1{2\varepsilon^2}
(1-|u|^2)^2,\quad E_{\varepsilon}(u)=\int_{\R^2}e_{\varepsilon}
(u)\,dx.
\endsplit
$$
Taking inner product of (2.11) with $\partial_t\,u_\varepsilon$,
one can obtain that for any $s>0$ in the maximal interval of
existence,
$$\int_{\R^2\times(0,s)} |\partial_t\,u_{\varepsilon}|^2\,dx\,dt + E_{\varepsilon}(u_{\varepsilon}(s))\leq
E(u_0).\tag 2.12
$$
Moreover,  repeating similar arguments in Lemma 3 (see below
(2.29) below) yields that the solution $u_{\varepsilon}$ belongs
to $\tilde{V}(0, T_{\varepsilon})$ for a maximum time
$T_{\varepsilon}$ and hence is regular $\R^2\times [0,
T_{\varepsilon})$. The maximum time $T_{\varepsilon}$ is
characterized in the following: For a singular point $x_0$ at
$T_{\varepsilon}$, there are $\varepsilon_0$ and $R_0>0$ such that
$$\limsup_{t\to T_{\varepsilon}} \int_{B_R(x_0)}|\nabla u_{\varepsilon}(\cdot ,t)|^2\,dx\geq
\varepsilon_0>0
$$
for any positive $R\leq R_0$.

Next, we will show that there is a uniform lower bound time
$t_1>0$ such that $T_{\varepsilon}\geq t_1$ and $u_{\varepsilon}$
is bounded in $\tilde{V}(0, t_1)$ uniformly in $\varepsilon$.

A similar argument as in Lemma 1 shows
$$\int_{B_R(x_0)}e_{\varepsilon} (u_\varepsilon (x,t))\,dx\leq \int_{B_{2R}(x_0)} e_{\varepsilon} (u_0 (x))\,dx
+C\frac t{R^2}\int_{\R^2}|\nabla u_0|^2\,dx
$$ for $t\leq T_{\varepsilon}$.

It follows from this inequality that for suitably small
$\varepsilon_1$ and $R_0$, there is a time $t_1$ uniform in
$\varepsilon$ with $t_1\leq T_{\varepsilon}$ such that
$$\sup_{0\leq t\leq t_1}\int_{B_R(x_0)}e_{\varepsilon} (u_\varepsilon (x,t))\,dx< \varepsilon_1\tag 2.13$$
for $R\leq R_0$ and thus $u_\varepsilon$ is smooth for $[0,t_1]$
for all $\varepsilon>0$. Next, we claim that for $0\leq t\leq t_1$
$$\frac 12 \leq|u_{\varepsilon}(x,t)|\leq \frac 32 \ \text{for all} \
x\in\R^2.
$$

To verify this claim, we re-scale the solution by $\tilde
u(x,t)=u_{\varepsilon} (\varepsilon x, \varepsilon^2 t)$. Then
$\tilde u$ satisfies
$$ \frac {\partial \tilde u^i}{\partial
t}= \nabla_{\a} \left [W_{p_{\a}^i} (\tilde u, \nabla \tilde u)
\right ]-W_{u^i}(\tilde u,\nabla \tilde u)+ \tilde{u}^i (1-|\tilde
u|^2)\tag 2.14
$$
with initial value $u_0(\varepsilon x)$. Let $\tau$ be the maximal
time in $[0,\frac {t_1}{\varepsilon^2}]$ such that (2.14) holds,
i.e.,
$$\frac{1}{2}\leq|\tilde u(x,t)|\leq \frac{3}{2}\tag 2.15$$
for any $(x,t)\in\R^2\times [0,\tau]$. Note that in this case, the
basic energy inequality (2.12) becomes
$$\split
& \int_{\R^2\times[0,s]}|\partial_t\,\tilde{u}|^2\,dx\,dt +
\int_{\R^2}(W(\tilde{u},\nabla\tilde{u})(s) +
\frac{1}{2}(1-|\tilde{u}(s)|^2)^2)dx \leq \ E(u_0)\\
& \text{\ for all} \ s\in\left[
0,\frac{t_1}{\varepsilon^2}\right],
\endsplit\tag 2.16
$$
and the condition (2.13) turns into
$$\esssup_{0\leq
s\leq\frac{t_1}{\varepsilon^2},x\in\R^2}
\int_{B_{\frac{R}{\varepsilon}}(x)} \left(
|\nabla\tilde{u}(\cdot,s)|^2+\frac{1}{2}(1-|\tilde{u}|^2)^2
\right)dx <\varepsilon_1.\tag 2.17
$$ for $R\leq R_0$.

Multiplying (2.14) by $\laplacian \tilde u$ and integrating over
$\R^2$ lead to
$$\split
\frac 1 2 &\frac d {dt}\int_{\R^2} |\nabla
 \tilde u|^2\,dx+\int_{\R^2} \nabla_{\a}  \left
[W_{p_{\a}^i} ( \tilde u,\nabla
 \tilde u) \right ]\laplacian  \tilde u^i\,dx\\
& -\int_{\R^2}W_{ \tilde u^i}( \tilde u,\nabla
 \tilde u)\laplacian  \tilde u^i \,dx +  \int_{\R^2}
 \tilde u^i (1-| \tilde u|^2)
\laplacian  \tilde u^i \,dx=0.
\endsplit\tag 2.18
$$
Note that
$$\split
 \int_{\R^2}
  \tilde u^i (1-| \tilde u|^2)
\laplacian  \tilde u^i \,dx &=- \int_{\R^2}|\nabla
 \tilde u|^2 (1-|\tilde u|^2)\,dx +2 \int_{\R^2} |\nabla | \tilde
 u|^2|^2\,dx.
\endsplit
$$
Then, combining the above identity with (2.18) yields that for any
$s,t\in [0,\tau]$ with $s\leq t$,
$$\split
&\quad \frac 1 2 \int_{\R^2} |\nabla \tilde u (\cdot
,t)|^2\,dx+\int_s^t \int_{\R^2} W_{p_{\a}^i p_{\g}^j} (\tilde
u,\nabla \tilde u)\nabla^2_{\a\b}\tilde u^i \nabla^2_{\g\b}\tilde
u^j \,dx\,dt
\\& \leq C\int_{\R^2} |\nabla u (\cdot
,s)|^2\,dx+ C\int_s^{t}\int_{\R^2}|\nabla \tilde u|^4+\eta
[(1-|\tilde u|^2)^2 +|\nabla^2 \tilde u|^2]\,dx\,dt
\endsplit\tag 2.19
$$
for a sufficiently small $\eta >0$ to be chosen.

On the other hand, it follows from (2.14) and (2.15) that
$$\int_{\R^2\times[s,t]} (1-|\tilde{u}|^2)^2\,dx\,dt\leq
C\int_{\R^2\times[s,t]} (|\partial_t\,\tilde{u}|^2 +
|\nabla\tilde{u}|^4 + |\nabla^2\tilde{u}|)dx\,dt.\tag 2.20
$$

Combining Lemma 2 with (2.17) shows that
$$\int_{\R^2\times[s,t]} |\nabla\tilde{u}|^4\,dx\,dt \leq C_1\,\varepsilon_1 \left(
\int_{\R^2\times[s,t]} |\nabla^2\tilde{u}|^2\,dx\,dt +
\frac{\varepsilon^2}{R^2_0} \int_{\R^2\times[s,t]}
|\nabla\tilde{u}|^2\,dx\,dt \right).\tag 2.21$$

As a consequence of (2.19)-(2.21), (2.16), and suitable choices of
$\eta$ and $\varepsilon_1$, one can get that
$\tilde{u}\in\tilde{V}_\tau$, and for any $0\leq s\leq t\leq\tau$,
$$\int_{\R^2\times[s,t]} \left[
|\nabla^2\tilde{u}|+(1-|\tilde{u}|^2)^2\right] dx\,dt \leq
CE(u_0)(1+\varepsilon^2(t-s)R^{-2}_0),\tag 2.22
$$
$$\int_{\R^2\times[s,t]} |\nabla\tilde{u}(x,t)|^4\,dx\,dt \leq
C\varepsilon_1\,E(u_0)(1+\varepsilon^2(t-s)R^{-2}_0).\tag 2.23
$$

By a similar argument as in the proof of Lemma 4, one can derive
from (2.14) that there exists a positive uniform constant $a$ such
that
$$\split
& \frac{1}{2}\frac{d}{dt}
\int_{\R^2}|\partial_t\,\tilde{u}(x,t)|^2\,dx +
a\int_{\R^2}|\partial_t\,\nabla\tilde{u}(x,t)|^2\,dx +
\frac{1}{2}\int_{\R^2}|\partial_t (|\tilde{u}(x,t)|^2)|^2\,dx\\
\leq &
C\int_{\R^2}|\partial_t\,\tilde{u}(x,t)|^2\,|\nabla\tilde{u}(x,t)|^2dx
+
\int_{\R^2}|\partial_t\,\tilde{u}(x,t)|^2\,(1-|\tilde{u}(x,t)|^2)dx
\endsplit\tag 2.24
$$

Note that
$$\split
& \
C\int_{\R^2}|\partial_t\,\tilde{u}(x,t)|^2\,|\nabla\tilde{u}(x,t)|^2 dx\\
\leq & \
C\left(\int_{\R^2}|\partial_t\,\tilde{u}(x,t)|^4\,dx\right)^{\frac{1}{2}}
\left(\int_{\R^2}|\nabla\tilde{u}(x,t)|^4\,dx\right)^{\frac{1}{2}}\\
\leq & \
C\left(\int_{\R^2}|\partial_t\,\tilde{u}(x,t)|^2\,dx\right)^{\frac{1}{2}}
\left(\int_{\R^2}|\partial_t\,\nabla\tilde{u}(x,t)|^2\,dx\right)^{\frac{1}{2}}
\left(\int_{\R^2}|\nabla\tilde{u}(x,t)|^4\,dx\right)^{\frac{1}{2}}\\
\leq & \
\frac{a}{2}\int_{\R^2}|\partial_t\,\nabla\tilde{u}(x,t)|^2\,dx +
C\left(\int_{\R^2}|\nabla\tilde{u}(x,t)|^4\,dx\right)
\int_{\R^2}|\partial_t\,\tilde{u}(x,t)|^2\,dx.
\endsplit
$$

Hence,
$$\split
& \frac{d}{dt}\int_{\R^2}|\partial_t\,\tilde{u}(x,t)|^2\,dx +
a\int_{\R^2}|\partial_t\,\nabla\tilde{u}(x,t)|^2\,dx +
\int_{\R^2}|\partial_t(|\tilde{u}(x,t)|^2)|^2\,dx\\
\leq & \left(C\int_{\R^2}|\nabla\tilde{u}(x,t)|^4\,dx\right)
\int_{\R^2}|\partial_t\,\tilde{u}(x,t)|^2\,dx +
2\int_{\R^2}(\partial_t\,\tilde{u}(x,t)|^2
(1-|\tilde{u}(x,t)|^2)dx,
\endsplit
$$
which yields immediately that for any $0\leq
t\leq\tau\in(0,\frac{t_1}{\varepsilon^2}]$,
$$\split
& \ \int_{\R^2}|\partial_t\,\tilde{u}(x,t)|^2\,dx\\
\leq & \
e^{C\int^t_0\!\!\int_{\R^2}|\nabla\tilde{u}(x,t)|^4\,dx\,dt}
\int_{\R^2}|\partial_t\,\tilde{u}(x,0)|^2\,dx\\
& \ + \int^t_0
\left(e^{C\int^t_0\!\!\int_{\R^2}|\nabla\tilde{u}(x,l)|^4\,dx\,dl}
\int_{\R^2}|\partial_t\,\tilde{u}(x,s)|^2 (1-|\tilde{u}(x,s)|^2)dx
\right)ds.
\endsplit
$$

It follows from this, (2.23), $u_0\in H^2_b$, and (2.14) that
$$\int_{\R^2}|\partial_t\,\tilde{u}(x,t)|^2\,dx \leq C_1\equiv
C_1(E(u_0),\varepsilon_1,t_1,||u_0||_{H^2_b}, R_0)\tag 2.25
$$
with a positive constant $C_1$ independent
$\tau\in(0,\frac{t_1}{\varepsilon^2})$ given by
$$C_1=C(||u_0||_{H^2_b}) \, e^{\varepsilon_1\,E(u_0)
(1+\frac{t_1}{R^2_0})}.\tag 2.26$$

Using (2.18), an integration by parts yields implies that for all
$t\in(0,\tau]$,
$$\int_{\R^2}|\nabla^2\tilde{u}(\cdot,t)|^2\,dx\leq
C\int_{\R^2}|\nabla\tilde{u}(x,t)|^4\,dx +
C\int_{\R^2}(1-|\tilde{u}|^2)^2\,dx +
C\int_{\R^2}|\partial_t\,\tilde{u}(x,t)|^2\,dx$$

Due to Lemma 2, and (2.17), one has
$$C\int_{\R^2}|\nabla\tilde{u}(x,t)|^4\,dx \leq
\varepsilon_1\int_{\R^2}|\nabla^2\,\tilde{u}(x,t)|^2\,dx +
\frac{C\varepsilon_1\,\varepsilon^2}{R^2_0} E(u_0).$$

Thus one can get that for all $t\in(0,\tau)$.
$$\int_{\R^2}|\nabla^2\,u(x,t)|^2\,dx \leq CE(u_0)
(1+\frac{\varepsilon_1\,\varepsilon^2}{R^2_0}) + CC_1.\tag 2.27
$$

By the Sobolev embedding theorem, $\tilde{u}$ is
$\beta$-H\"{o}lder continuous in $x$ uniformly in all
$t\in[0,\tau]$ with $\beta<1$. Repeating the similar analysis as
in the proof of Lemma 4 and using Proposition 13 in Appendix, we
get $\tilde{u}\in C^{1,\frac{1}{8}}$ on $\R^2\times(0,\tau)$. If
there is a $x_1\in\R^2$ such that either
$|\tilde{u}(x_1,t)|<\frac{1}{2}$ or
$|\tilde{u}(x_1,t)|>\frac{3}{2}$. By the uniform H\"{o}lder
continuity of $\tilde{u}$, there exists a constant $C_2$ with the
property that $\frac{1}{4C_2}<\frac{R_0}{\varepsilon}$, and
$$(1-|\tilde{u}(x,t)|^2)^2 \geq \frac{1}{4}, \qquad x\in
B_{\frac{1}{4C_2}}(x_1).$$

Hence,
$$\int_{B_{\frac{1}{4C_2}}(x_1)} (1-|\tilde{u}(x,t)|^2)^2\,dx \geq
\frac{1}{4}|B_{\frac{1}{4C_2}}(0)|>2\varepsilon_1.\tag 2.28
$$
which contradicts to (2.17) for suitably small $\varepsilon_1$.
Here we have used the fact that $C_2$ depends only the upper bound
of $C_1$, which may be chosen to be independent of $\varepsilon_1$
by the choice of $t_1$. This implies that $\frac 1 2 \leq |\tilde
u(x,t)|\leq \frac 32$ for all $t\in [0,\tau ]$. By the continuity
of $u$ at $\tau$ and the maximal choice of $\tau$, $\tau$ must be
the value $\frac {t_1}{\varepsilon^2}$. This shows that (2.14)
holds for all $t\in [0,t_1]$.

Next, it follows from (2.12) and (2.22)-(2.23) that
$u_{\varepsilon}$ are uniformly bounded in $\tilde{V}(0,t_1)$ for
all $\varepsilon$ and

$$\int_0^{t_1}\!\!\int_{\R^2}\left[ |\nabla^2\,u_{\varepsilon}(x,t)|^2+ \frac{1}{\varepsilon^2}
(1-|u_\varepsilon(x,t)|^2)^2\right]dx\,dt \leq CE(u_0)
(1+\frac{t_1}{R^2_0}),\tag 2.29 $$
$$
\int^{t_1}_0\!\!\int_{\R^2} |\nabla\,u_\varepsilon(x,t)|^4\,dx\,dt
\leq C\varepsilon_1\,E(u_0) (1+\frac{t_1}{R^2_0}).\tag 2.30
$$

Letting $\varepsilon\to 0$, we can prove local existence of a
solution of (1.4) in $V(0,t_1)$.\qed
\enddemo

Now we complete the proof of Theorem A.

\demo{Proof of Theorem A} By Theorem 7, there is a local solution
$u$ on $[0, t_1)$ for some $t_1>0$. By Lemma 3 and Lemma 4, the
solution can be extended to $[0,T_1)$ for a maximal time $T_1>0$
such that there is a singular set $\Sigma$ at $T_1$. Each
singularity $ x_{i}^{1} \in\Sigma$ at $T_1$ is characterized by
the condition
$$\limsup_{t \nearrow  T_{1}} E\left( u\left(
t\right)  ,B_{R}\left( x_{i}^{1}\right) \right) \geq
\varepsilon_{0}$$ for any $R>0$ with $R\leq R_0$. It is easy to
see the solution $u\in V$ is  regular for all $t\in (0, T_1)$. By
Lemma 1, we can show that the singular set $\Sigma$ and the
singular times are finite (See \cite {St1}). Theorem A is thus
proved. \qed
\enddemo

\proclaim{Remark} There is an open problem to prove the uniqueness
of the weak solutions. But, we can prove the uniqueness of smooth
solutions (see Lemma 11 below in Section 3).
\endproclaim

\head {\bf 3. Global existence for  the Ericksen-Leslie
system}\endhead

In this section, we derive a-priori estimates  for solutions to
the Ericksen-Leslie system (1.5)-(1.7). Without loss of
generality, we assume that $\nu =\lambda =1$ in (1.5).

For the case $\Omega =\R^2$, we still consider (1.5)-(1.7) in
$\R^3$ by taking $\frac {\partial v}{\partial x_3}=0, \frac
{\partial u}{\partial x_3}=0$. In this case, $\nabla \cdot v =
\frac {\partial v^1}{\partial x_1}+ \frac {\partial v^2}{\partial
x_2}=0$ in (1.6) is well-defined.

 For two positive constants $\tau$ and
$T$ with $\tau < T$, we denote
 $$\split
 V(\tau, T):= \{ u: \R^2\times &[\tau , T]\to S^2\, |\quad  u \quad \text
  {is measureable and satisfies}\\
&
 \esssup_{\tau\leq t\leq T} \int_{\R^2} |\nabla u(\cdot ,t)|^2\, dx
 +\int_{\tau}^T \int_{\R^2} |\nabla^2 u|^2 +|\partial_t u|^2 \,dx\,dt <\infty
  \}
 \endsplit
 $$
 and
 $$
 \split H(\tau , T):= \{ v: \R^2\times &[\tau , T]\to \R^2\, |\quad  v \quad \text
  {is measureable and satisfies}\\
&
 \esssup_{\tau\leq t\leq T} \int_{\R^2} |v(\cdot ,t)|^2\, dx
 +\int_{\tau}^T \int_{\R^2} |\nabla  v|^2 \,dx\,dt <\infty
  \}.
 \endsplit
 $$
 For each pair $(u,v)$, define
$$e(u,v)=W(u,\nabla u)+\frac 12 |v|^2, \qquad
E(u,v)=\int_{\R^2}e(u,v)dx.
$$

\proclaim {Lemma 7} Let $(u, v)\in V(0 ,T)\times H(0,T)$ be a
solution of (1.5)-(1.7) with initial  values   $u_0\in H^1(\R^2;
S^2)$ and $v_0\in L^2(\R^2; \R^3)$. Then for $t\in(0,T]$,
$$ \split & \int_{\R^2} e( u (\cdot, t), v (\cdot, t)) \,dx +
 \int_0^t\int_{\R^2} (| u_t +(v\cdot \nabla ) u|^2+ |\nabla v|^2) \,dx\,dt\\
&=\int_{\R^2} e(u_0, v_0) \,dx.\endsplit \tag 3.1
$$
\endproclaim
\demo {Proof} Multiplying (1.5) by $v$ and using (1.6), one gets
$$\frac 1 2\frac d{dt}\int_{\R^2} |v|^2\,dx +\int_{\R^2} |\nabla v|^2\,dx =  \int_{\R^2}
\nabla_{j}v^i\nabla_{i}u^kW_{p_j^k}(u, \nabla u)\,dx.\tag 3.2
$$
Multiplying (1.7) by $u_t +(v\cdot \nabla )u $ yields
$$\split
&\quad \int_{\R^2}(u^i_t +(v\cdot\nabla )u^i) \big ( \nabla_{\a}
\left [W_{p_{\a}^i} (u,\nabla u)- u^ku^i V_{p_{\a}^k}(u,\nabla u)
\right
] )\,dx \\
& +\int_{\R^2}(u^i_t +(v\cdot\nabla )u^i) (-W_{u^i}(u,\nabla
u)+W_{u^k}(u,\nabla u)u^k u^i \big )\,dx\\
& +\int_{\R^2}(u^i_t +(v\cdot\nabla )u^i) \big
(+W_{p_{\a}^k}(u,\nabla u)\nabla_{\a}u^k u^i
+V_{p_{\a}^k}(u,\nabla u) u^k\nabla_{\a}u^i \big )\,dx\\
=&\int_{\R^2} \left |u_t +(v\cdot \nabla )u \right |^2\,dx .
\endsplit\tag 3.3
$$
Note that $|u|^2=1$ implies
$$u^i\partial_t u^i=0,\quad u^i\nabla_{x_{\a}} u^i=0. $$
Integration by parts
yields
$$\split
\int_{\R^2}u^i_t \big (&\nabla_{\a} \left [W_{p_{\a}^i} (u,\nabla
u)- u^ku^i V_{p_{\a}^k}(u,\nabla u) \right ]-W_{u^i}(u,\nabla
u)+V_{p_{\a}^k}(u,\nabla u) u^k\nabla_{\a}u^i
\big )\,dx\\
&= -\int_{\R^2} \nabla_{\a}u^i_t  \left [W_{p_{\a}^i} (u,\nabla
u)\right ]-u^i_t  W_{u^i}(u,\nabla u)\,dx=- \frac d{dt}
\int_{\R^2} W(u,\nabla u)\,dx.
\endsplit\tag 3.4
$$
Using (1.6) and integrating by parts, we  get
$$\split
\int_{\R^2} (v\cdot \nabla )u^i ( &\nabla_{\a} \left [W_{p_{\a}^i}
(u,\nabla u)- u^ku^i V_{p_{\a}^k}(u,\nabla u) \right
]-W_{u^i}(u,\nabla u)\\
&\quad +V_{p_{\a}^k}(u,\nabla u) u^k\nabla_{\a}u^i  )\,dx\\
=-\int_{\R^2}\nabla_{\a} v^k\nabla_{k} u^i &W_{p_{\a}^i} (u,\nabla
u)+v^k[\nabla_k\nabla_{\a}u^i W_{p_{\a}^i}(u,\nabla u)
+ \nabla_{k} u^iW_{u^i}(u,\nabla u) ]\,dx\\
=- \int_{\R^2}\nabla_{\a} v^k\nabla_{k} u^i &W_{p_{\a}^i}
(u,\nabla u) \,dx.
\endsplit\tag 3.5
$$
It follows from (3.3)-(3.5) that
$$\split
\frac d {dt}\int_{\R^2} W(u,\nabla u)\,dx + \int_{\R^2} | u_t
+(v\cdot \nabla ) u|^2\,dx  = - \int_{\R^2}\nabla_{\a}
v^k\nabla_{k} u^i &W_{p_{\a}^i} (u,\nabla u)\,dx.\endsplit \tag
3.6
$$
Therefore, (3.1) follows from integrating (3.2) and (3.6) in $t$.
\qed
\enddemo

By the same proof as in \cite {St1; Lemma 3.1}, there exists a
constant $C_1$ such that for any $f\in H(0 ,T)$ and any $R>0$, it
holds that
$$\split  \quad \int_{\R^2\times [0 ,T]} |f|^4\,dx\,dt\leq &C_1\esssup_{0\leq t\leq T, x\in \R^2}
\int_{B_R(x)}|f(\cdot ,t)|^2\,dx  \\& \cdot \left
(\int_{\R^2\times [0 ,T]}
 |\nabla f|^2  \,dx\,dt +R^{-2} \int_{\R^2\times [0 ,T]} |f|^2
\,dx\,dt \right ).
\endsplit\tag 3.7
$$

Then, we have
 \proclaim {Lemma 8} Let $(u, v)\in V(0 ,T)\times
H(0,T)$ be a solution of (1.5)-(1.7) with initial  values $u_0\in
H^1$ and $v_0\in L^2$. Then there are constants $\varepsilon_1$
and $R_0>0$ such that if
$$\esssup_{0\leq t\leq T, x\in \R^2} \int_{B_R(x)} e( u(\cdot ,t),v(\cdot ,t))\,dx<\varepsilon_1
 $$ for  any
$R\in (0, R_0]$, then
$$\int_{\R^2\times [0 ,T]} |\nabla^2 u|^2+ |\nabla v|^2\,dx\,dt\leq C
  \, (1 +T R^{-2})\int_{\R^2} e(u_0,v_0)\,dx,\tag 3.8
$$
$$\int_{\R^2\times[0,T]}(|\nabla u|^4 +|v|^4)dx\,dt \leq
C\varepsilon_1 (1+TR^{-2}) \int_{\R^2} e(u_0,v_0)dx\tag 3.9$$
\endproclaim
\demo {Proof}Multiplying $\laplacian u^i$ with (1.7) yields
$$\split
&\quad \int_{\R^2}(\frac {\partial u^i}{\partial t}+ (v\cdot
\nabla ) u^i )\laplacian u^i\,dx\\
= & \int_{\R^2}  \nabla_{\a} \left [W_{p_{\a}^i} (u,\nabla
u)- u^ku^i V_{p_{\a}^k}(u,\nabla u) \right ]\laplacian u^i\,dx\\
& -\int_{\R^2}W_{u^i}(u,\nabla u)\laplacian u^i dx + \int_{\R^2}
u^k\,u^i\,W_{u^k}(u,\nabla u)\laplacian u^i\,dx
\\
& +\int_{\R^2}W_{p_{\a}^k}(u,\nabla u)\nabla_{\a}u^k u^i
\laplacian u^i \,dx +\int_{\R^2}V_{p_{\a}^k}(u,\nabla u)
u^k\nabla_{\a}u^i \laplacian u^i\,dx.
\endsplit
$$

As in the proof of Lemma 3, one can derive
$$\split
&\quad \frac {d}{dt}\int_{\R^2} |\nabla u|^2\,dx +\int_{\R^2}
W_{p_{\a}^i p_{\g}^j} (u,\nabla u)\nabla^2_{\a\b}u^i
\nabla^2_{\g\b}u^j \,dx\\
&\leq C\int_{\R^2} (|\nabla u|^2+|v|^2) (|\nabla
u|^2+|\nabla^2u|)\,dx \\
&\leq \frac b 4 \int_{\R^2}|\nabla^2 u|^2\,dx + C \int_{\R^2}
(|\nabla u|^4+|v|^4) \,dx . \endsplit
$$

Applying (3.7) and Lemma 2 again shows
$$\split & \quad \int_{\R^2\times [0 ,T]} |\nabla u|^4+|v|^4\,dx\,dt\\
\leq & C_1\varepsilon_1 \int_{\R^2\times [0 ,T]}
 |\nabla^2 u|^2  +|\nabla v|^2 \,dx\,dt + C_1\,\varepsilon_1\,R^{-2} \int_{\R^2\times [0 ,T]} |\nabla
 u|^2+|v|^2
\,dx\,dt.
\endsplit
$$
Then (3.8) and (3.9) follow by choosing $\varepsilon_1=\frac
b{4C_1}$.\qed\enddemo

\proclaim {Lemma 9} Let $(u,v)$ be a solution of (1.5)-(1.7) with
initial values $(u_0,v_0)$ with $u\in V(0 ,T)$ and $v\in H(0,T)$.
Assume that there exist  constants $\varepsilon_1>0$ and $R_0>0$
such that
 $$\sup_{ x\in \R^2, 0\leq t\leq T} \int_{B_{R_0}(x)} |\nabla u (x,t)|^2 +|v\left(\cdot , t\right)|^2 \,dx
 <\varepsilon_1.$$
Then for all $t\in [0,T]$, $x_0\in \R$ and $R\leq R_0$, it holds
that
$$\split
&\int_{B_R(x_0)} e(u (\cdot,t),v(\cdot,t))\,dx +
\int^t_0\!\!\int_{B_R(x_0)} (|\nabla v|^2 +
\frac{1}{2}|\partial_t\,u +v\cdot\nabla u|^2) dx\,dt\\
\leq & \int_{B_{2R}(x_0)} e(u_0,v_0)\,dx +C_2
\frac{t^{\frac{1}{2}}}{R}(1+\frac{t}{R^2})^{\frac{1}{2}}
\int_{\R^2} e(u_0,v_0)dx,
\endsplit\tag 3.10
$$
where  $C_2$ is a uniform positive constant.
\endproclaim
\demo{Proof}   Let $\phi\in C_0^{\infty} (B_{2R} (x_0))$ be a
cut-off function with $\phi\equiv 1$ on $B_R(x_0)$ and $|\nabla
\phi |\leq \frac C R$, $|\nabla^2 \phi|\leq \frac C{R^2}$ for all
$R\leq R_0$.

Multiplying (1.5) by $v\phi^2$ and integrating show
$$\split
&\quad \int_{\R^2} v_t\cdot v\phi^2 + (v\cdot \nabla )v \cdot
v\phi^2
-\laplacian v \cdot v\phi^2 +\nabla P\cdot v \phi^2 \,dx\\
&=\int_{\R^2}  \nabla_{x_i} u^k W_{p_{j}^k} (u,\nabla
u)\nabla_{x_j} v^i\phi^2\,dx +\int_{\R^2}  \nabla_{x_i} u^k
W_{p_{j}^k} (u,\nabla u) v^i\nabla_{x_j}\phi^2\,dx.
\endsplit
$$
Integrating by parts yields
$$\split
&\int_{\R^2} (v_t\cdot v\phi^2-\laplacian v \cdot
v\phi^2)\,dx=\frac 1 2 \frac d{dt}\int_{\R^2} |v|^2\phi^2
\,dx+\int_{\R^2} \nabla v \cdot \nabla (v\phi^2)\,dx\\
&=\frac 1 2 \frac d{dt}\int_{\R^2} |v|^2\phi^2 \,dx+\int_{\R^2}
|\nabla v|^2 \phi^2 \,dx -\int_{\R^2}|v|^2(|\nabla \phi|^2
+\phi\laplacian \phi) \,dx.
\endsplit
$$
Integrating by parts and using (1.6) give
$$\int_{\R^2} \nabla_{x_i} P v^i \phi^2\,dx =-2\int_{\R^2} P v^i\phi\nabla_{x_i}
\phi\,dx$$ and
$$\int_{\R^2} v^k \nabla_{x_k} v^i v^i\phi^2=\frac 12 \int_{\R^2} v^k
\nabla_{x_k} (|v|^2) \phi^2=-\int_{\R^2} v^k |v|^2 \phi
\nabla_{x_k}\phi \,dx.
$$

Hence,
$$\split
&\quad \frac 1 2\frac d {dt}\int_{\R^2} |v|^2\phi^2\,dx  + \int
|\nabla v|^2
\phi^2\,dx\\
&  =\int_{\R^2} ( |v|^2+2P+|\nabla u|^2) v^i\phi \nabla_{x_i}\phi
\,dx  + \int_{\R^2}|v|^2(|\nabla \phi|^2 +\phi\laplacian \phi)
\,dx\\
&+\int_{\R^2}  \nabla_{x_i} u^k W_{p_{j}^k} (u,\nabla
u)\nabla_{x_j} v^i\phi^2\,dx +\int_{\R^2}  \nabla_{x_i} u^k
W_{p_{j}^k} (u,\nabla u) v^i\nabla_{x_j}\phi^2\,dx.
\endsplit\tag 3.11
$$
Multiplying (1.7) by $(u_t^i +(v\cdot \nabla )u^i )\phi^2$  and
using $|u|=1$ lead to
$$\split
&\quad \int_{\R^2} |u_t +(v\cdot \nabla )u|^2\phi^2\,dx\\
= & \int_{\R^2}  (u^i_t+v^l\nabla_l u^i ) \nabla_{\a}
[W_{p_{\a}^i} (u,\nabla u)-u^ku^i V_{p_{\a}^k}(u,\nabla
u)]\varphi^2\,dx\\
&+ \int_{\R^2} (u^i_t+v^l\nabla_l u^i )(-W_{u^i}(u,\nabla
u)+V_{p_{\a}^k }(u,\nabla u) u^k \nabla_{\a}u^i )\phi^2\,dx.
\endsplit
$$
Integration by parts yields
$$\split
& \int_{\R^2} u^i_t\nabla_{\a} [W_{p_{\a}^i} (u,\nabla u)-u^ku^i
V_{p_{\a}^k}(u,\nabla u)]\phi^2\,dx\\
&+\int_{\R^2} u^i_t(-W_{u^i}(u,\nabla
u)+V_{p_{\a}^k }(u,\nabla u) u^k \nabla_{\a}u^i )\phi^2\,dx\\
&=-\int_{\R^2} [\nabla_{\a} u_t^i W_{p_{\a}^i} (u,\nabla u)+u_t^i
W_{u^i}(u,\nabla u)]\phi^2\,dx-\int_{\R^2} u_t^iW_{p_{\a}^i} (u,\nabla u)\nabla_{\a}\phi^2 \,dx\\
&=-\frac d{dt}\int_{\R^2}W(u,\nabla u)\phi^2\,dx-2\int_{\R^2}
u_t^iW_{p_{\a}^i} (u,\nabla u)\phi \nabla_{\a}\phi \,dx.
\endsplit
$$
Integrating by parts twice and using (1.6), we obtain
$$\split
& \int_{\R^2} v^l\nabla_l u^i \nabla_{\a} [W_{p_{\a}^i} (u,\nabla
u)-u^ku^i V_{p_{\a}^k}(u,\nabla u)]\phi^2\,dx\\
&\quad +\int_{\R^2} v^l\nabla_l u^i(-W_{u^i}(u,\nabla
u)+V_{p_{\a}^k }(u,\nabla u) u^k \nabla_{\a}u^i )\phi^2\,dx\\
&=- \int_{\R^2} (\nabla_{\a} v^l\nabla_l u^i +v^l\nabla_l
\nabla_{\a}u^i) W_{p_{\a}^i} (u,\nabla u)\phi^2\,dx\\
&-\int_{\R^2} ( v^l\nabla_l u^i) W_{p_{\a}^i} (u,\nabla
u)\nabla_{\a} \phi^2\,dx -\int_{\R^2} v^l\nabla_l
u^iW_{u^i}(u,\nabla
u)\phi^2\,dx\\
&=- \int_{\R^2} \nabla_{\a} v^l\nabla_l u^i W_{p_{\a}^i} (u,\nabla
u)\phi^2\,dx-\int_{\R^2} ( v^l\nabla_l u^i) W_{p_{\a}^i} (u,\nabla
u)\nabla_{\a} \phi^2\,dx\\
&\quad +2\int_{\R^2} v^l W(u,\nabla u)\phi \nabla_l\phi\,dx.
\endsplit
$$
Combing above three identities yields
$$\split
&  \frac d {dt} \int_{\R^2} W(u,\nabla u) \phi^2\,dx+\int_{\R^2}
| u_t +(v\cdot \nabla ) u|^2
\phi^2\,dx\\
& =- \int_{\R^2} \nabla_{\a} v^k\nabla_k u^i W_{p_{\a}^i}
(u,\nabla u)\phi^2\,dx-\int_{\R^2} (u_t^i+ v^k\nabla_k u^i)
W_{p_{\a}^i} (u,\nabla
u)\nabla_{\a} \phi^2\,dx\\
&\quad +2\int_{\R^2} v^k W(u,\nabla u)\phi \nabla_k\phi\,dx\\
&\leq - \int_{\R^2} \nabla_{\a} v^k\nabla_k u^i W_{p_{\a}^i}
(u,\nabla u)\phi^2\,dx +\frac 1 2 \int_{\R^2} | u_t +(v\cdot
\nabla ) u|^2 \phi^2\,dx\\
&\quad +C\int_{\R^2} |\nabla u|^2 |\nabla \phi |^2\,dx +
2\int_{\R^2} v^k\,W(u,\nabla u)\phi\,\nabla_k\,\phi\,dx.
\endsplit\tag 3.12
$$
Integrating (3.11) and (3.12) in $t$ on $[0,s]$ leads to
$$\split
&\int_{\R^2} e(u (\cdot ,s), v (\cdot ,s)) \phi^2 \,dx
+\int_0^s\int_{\R^2} (|\nabla v|^2 + \frac 1 2 | u_t +(v\cdot
\nabla ) u|^2 )\phi^2\,dxdt
\\
&\leq  \int_{\R^2} e(u_0,v_0)\phi^2 \,dx
+\int_0^s\int_{\R^2}(|v|^2 +|\nabla u|^2+ 2P )v^i\phi \nabla_{x_i}\phi \,dxdt\\
& + 2\int_0^{s}\int_{\R^2} (v\cdot \nabla )u^k W_{p^k_i}(u,\nabla
u) \phi \nabla_{x_i}\phi \,dxdt +2\int_0^s\int_{\R^2} v^l
W(u,\nabla
u)\phi \nabla_l\phi\,dx\\
&+C\int_0^s \int_{\R^2} (|v|^2+|\nabla u|^2) (|\nabla \phi|^2
+|\phi |\,|\laplacian \phi |) \,dxdt.
\endsplit\tag 3.13$$

This, together with (3.1), shows immediately that
$$\split
&  \int_{B_R(x_0)} (|v (\cdot ,s)|^2 + |\nabla u (\cdot ,s)|^2)
\,dx +\int_0^s\int_{\R^2} (|\nabla v|^2 +\frac 1 2 |u_t +(v\cdot
\nabla ) u|^2 )\phi^2\,dxdt
\\
&\leq  \int_{B_{2R}(x_0)} (|v_0|^2 +|\nabla u_0|^2) \,dx
+C\int_0^s\int_{\R^2}(|v|^2 +|\nabla u|^2+ |P| ) |v| |\phi|  |\nabla \phi | \,dxdt\\
&  \quad   + C \frac s{R^2}\int_{\R^2} (|v_0|^2 +|\nabla u_0|^2)
\,dx.
\endsplit\tag 3.14
$$

It follows from H\"{o}lder inequality, (3.1) and (3.9) that
$$\split
& \ \int_0^s\!\!\int_{\R^2}(|v|^2 +|\nabla u|^2) |v| |\phi|
|\nabla \phi | \,dxdt\\
\leq & \ C\int^s_0\!\!\int_{\R^2}(|v|^2+|\nabla
u|^2)\frac{|v|}{R}dx\,dt\\
\leq & \ C \left( \int_0^s\!\!\int_{\R^2}(|v|^4 +|\nabla
u|^4)dx\,dt \right)^{\frac{1}{2}} \left( \int_0^s\!\!\int_{\R^2}
\frac{|v|^2}{R^2}dx\,dt \right)\\
\leq & \ C\varepsilon^{\frac{1}{2}}\frac{s^{\frac{1}{2}}}{R}
\left( 1+\frac{s}{R^2} \right)^{\frac{1}{2}}
\int_{\R^2}e(u_0,v_0)dx.
\endsplit \tag 3.15
$$

Similarly,
$$\split
\int_0^s\!\!\int_{\R^2}|P|  |v| |\phi|  |\nabla \phi | \,dxdt \leq
& C\int^s_0\!\!\int_{\R^2}|p|\frac{|v|}{R}dx\,dt\\
\leq & C \left( \int_0^s\!\!\int_{\R^2}|P|^2\,dx\,dt
\right)^{\frac{1}{2}}
\left( \int^s_0\!\!\int_{\R^2}\frac{|v|^2}{R^2}dx\,dt \right)^{\frac{1}{2}}\\
\leq & C\frac{s^{\frac{1}{2}}}{R} \left( \int_{\R^2} e(u_0,v_0)dx
\right)^{\frac{1}{2}} \left( \int^s_0\!\!\int_{\R^2}|P|^2 \,dxdt
\right)^{\frac{1}{2}}
\endsplit\tag 3.16
$$
for $R\leq R_0$.

On the other hand, it follows from the relation that
$$\laplacian   P=-\nabla_{x_ix_j} \left [ \nabla_{x_i}  u^k W_{p_j^k} (u,\nabla u)+  v^j v^i\right
]\quad \text {on } \R^2\times (0,T],
$$
due to (1.5), and the Calderon-Zygmund estimate (cf. \cite {CKN})
that
$$\int^s_0\!\!\int_{\R^2} |P|^2\,dx\,dt \leq C\int^s_0\!\!\int_{\R^2} (|\nabla
u|^4 +|v|^4)\,dx\,dt \leq C\,\varepsilon_1
\left(1+\frac{s}{R^2}\right) \int_{\R^2} e(u_0,v_0)dx.
$$
This, together with (3.16), yields
$$\int_0^s\!\!\int_{\R^2} |P|\,|v|\,|\varphi|\,|\nabla\phi|\,dx\,dt \leq C \frac{s^{\frac{1}{2}}}{R} \left( 1+\frac{s}{R^2} \right)^{\frac{1}{2}}
\, \varepsilon^{\frac{1}{2}}_1 \int_{\R^2} e(u_0,v_0) dx.\tag 3.17
$$
The desired estimate (3.10) now follows from (3.14), (3.15) and
(3.17). \qed
\enddemo

\proclaim {Lemma 10} Let $u\in V(0 ,T)$ and $v\in H(0,T)$ be a
solution of (1.5)-(1.7) with initial value $(u_0,v_0)\in
H_b^1(\R^2,S^2)\times L^2(\R^2,\R^3)$ and $\text {div }v_0 =0$.
Assume that there are constants $\varepsilon_1$ and $R_0>0$ such
that
$$\esssup_{0\leq t\leq T, x\in \R^2} \int_{B_R(x)}|\nabla
u(\cdot ,t)|^2+|v\left(\cdot , t\right)|^2 \,dx<\varepsilon_1
 $$ for  any
$R\in (0, R_0]$. Let $\tau$ be any positive constant. Then, for
$t\in [\tau, T]$, it holds that
$$\int_{\R^2} |\nabla^2 u(x,t) |^2+|\nabla v (x,t)|^2\,dx\leq C\,\tau^{-1} (1 +T
R^{-2}).\tag 3.18
$$
Moreover, $u$ and $v$ are regular for all $t\in (0,T)$.
\endproclaim
\demo {Proof} Note that, in a-priority, $\int_{\R^2} |\laplacian
v|^2$ and $\int_{\R^2}|\nabla^3 u|^2$ might not be finite.
However, by a standard cut-off argument, we can assume that
$\int_{\R^2} |\laplacian v|^2$ and $\int_{\R^2}|\nabla^3 u|^2$ are
finite without loss of generality in the following proof.

Multiplying (1.5) by $\laplacian v^i$ and integrating by parts, we
obtain
$$\split
&\quad \frac 1 2\frac d {dt} \int_{\R^2} |\nabla v|^2\,dx +\int_{\R^2} |\laplacian v|^2\,dx\\
& =\int_{\R^2} (v\cdot \nabla v^i)  \laplacian v^i\,dx
+\int_{\R^2} \nabla_{j}[\nabla_i u^k W_{p_j^k} (u,\nabla u) ] \laplacian v^i\,dx\\
&\leq \frac 1 4 \int_{\R^2} |\laplacian v|^2 \,dx+C\int_{\R^2}
|v\cdot\nabla v|^2\,dx + C\int_{\R^2}( |\nabla^2 u|^2+ |\nabla
u|^4) |\nabla u|^2\,dx .
\endsplit\tag 3.19
$$

Differentiating (1.7) in $x_{\b}$, multiplying  the above equation
by $\nabla_{\b} \laplacian u^i$ and then integrating by parts, one
can obtain
$$\split
& \quad -\frac  12 \frac d{dt}\int_{\R^2}|\laplacian u|^2\, dx+
\int_{\R^2}\left [ ( \nabla_{\b} v\cdot \nabla) u^i
 +(v \cdot \nabla )\nabla_{\b}u^i\right ] \nabla_{\b} \laplacian u^i\, dx\\
=& \int_{\R^2}\left [\nabla_{\b} \nabla_{\a} \left [W_{p_{\a}^i}
(u,\nabla u)- u^ku^i V_{p_{\a}^k}(u,\nabla u) \right ]-\nabla_{\b}
W_{u^i}(u,\nabla
u)\right ]\nabla_{\b} \laplacian u^i\,dx\\
&+ \int_{\R^2} \nabla_{\b} [W_{u^k}(u,\nabla u)u^k
u^i+W_{p_{\a}^k}(u,\nabla u)\nabla_{\a}u^k u^i ] \,\nabla_{\b}
\laplacian u^i\,dx\\
&+ \int_{\R^2} \nabla_{\b} [V_{p_{\a}^k}(u,\nabla u)
u^k\nabla_{\a}u^i] \,\nabla_{\b} \laplacian u^i\,dx  .
\endsplit\tag 3.20
$$
The first term on the righthand side of (3.20) is a bit more
complicated. Since $W(u,p)$ is quadratic in $p$, we have
$$\split
\nabla^2_{\g \b}W_{p_{\a}^i} (u,\nabla u)&= \nabla_{\g} [
W_{u^jp_{\a}^i}(u,\nabla u) \nabla_{\b}
u^j+W_{p_{\a}^i}(u,\nabla \nabla_{\b}u) ]\\
&=W_{u^jp_{\a}^i}(u,\nabla u) \nabla^2_{\g\b}
u^j+W_{u^ju^kp_{\a}^i}(u,\nabla u) \nabla_{\g} u^k\nabla_{\b}
u^j\\
&+W_{p_l^jp_{\a}^i}(u,\nabla \nabla_{\b}u) \nabla^3_{\b \g l}
u^j+W_{u^jp_{\a}^i}(u,\nabla \nabla_{\b}u) \nabla_{\g} u^j.
\endsplit \tag 3.21$$
Then, integrating by parts  and  using Young's inequality, we have
$$\split
&\int_{\R^2}\nabla_{\b} \nabla_{\a} W_{p_{\a}^i} (u,\nabla u)
\nabla_{\b} \laplacian u^i\,dx=\int_{\R^2}\nabla^2_{\g \b}
W_{p_{\a}^i} (u,\nabla u)
\nabla^3_{\g \b \a} u^i\,dx\\
& \geq \frac a 4\int_{\R^2} |\nabla^3 u|^2\,dx
-C\int_{\R^2}|\nabla u|^2(|\nabla u|^4+|\nabla^2u|^2 ) \,dx.
\endsplit\tag 3.22
$$
Note that $|u|^2=1$ implies
$$u^i\nabla_{\b} \laplacian u^i +\nabla_{\b} u^i \laplacian u^i =-\nabla_{\b} |\nabla u|^2. $$
By this identity, one can estimate the second term and the last
term on the righthand of (3.20) as follows:
$$\split
&\int_{\R^2}\nabla_{\b} [\nabla_{\a}\left ( u^ku^i
V_{p_{\a}^k}(u,\nabla u) \right ) -V_{p_{\a}^k}(u,\nabla u)
u^k\nabla_{\a}u^i] \nabla_{\b} \laplacian u^i\,dx\\
= & \int_{\R^2} \nabla^2_{\a\b}\left ( u^k V_{p_{\a}^k}(u,\nabla
u) \right ) u^i\, \nabla_{\b} \laplacian
u^i\,dx\\
&+\int_{\R^2}\nabla_{\b} u^i\,\nabla_{\a}\left ( u^k
V_{p_{\a}^k}(u,\nabla u) \right ) \nabla_{\b} \laplacian
u^i\,dx\\
\leq & \frac a 4 \int_{\R^2} |\nabla^3 u|^2\,dx +
C\int_{\R^2}|\nabla u|^2( |\nabla^2 u|^2+|\nabla u|^4) \,dx.
\endsplit\tag 3.23
$$
The other terms can be estimated easily in (3.20). Then it follows
from (3.20)-(3.23) that
$$\split
&\frac 1 2 \frac d{dt} \int_{\R^2} |\laplacian u|^2 dx + \frac a
4\int_{\R^2} |\nabla^3 u
|^2\,dx \\
\leq & C\int_{\R^2} |\nabla v|^2 |\nabla u|^2 + |v|^2 |\nabla^2
u|^2 +(|\nabla^2 u|^2+|\nabla u|^4) |\nabla u|^2 \,dx.
\endsplit\tag 3.24
$$
It follows from $-u\cdot \laplacian u=|\nabla u|^2$, (3.18) and
(3.24) that
$$\split
& \frac {d}{dt}\left ( \int_{\R^2} |\nabla^2 u|^2 +|\nabla v|^2
\right )+ \frac a 4 \int_{\R^2} \left (|\nabla^3 u|^2
 + |\nabla^2
v|^2 \right )\,dx\\
\leq & C\int_{\R^2} (|v|^2+|\nabla u|^2) \left ( |\nabla v|^2
+|\nabla^2 u|^2 \right )\,dx.
\endsplit\tag 3.25
$$
By the Gagliardo-Nirenberg-Sobolev inequality, one has
$$\split
& \ C\int_{\R^2} (|v|^2 + |\nabla u|^2) (|\nabla v|^2 + |\nabla^2
u|^2)dx\\
\leq & \ C\left( \int_{\R^2} (|v|^4 + |\nabla u|^4)dx
\right)^{\frac{1}{2}} \left( \int_{\R^2} (|\nabla v|^4 + |\nabla^2
u|^4)dx \right)^{\frac{1}{2}}\\
\leq & \ C\left( \int_{\R^2} (|\nabla v|^2 + |\nabla^2 u|^2)dx
\right)^{\frac{1}{2}} \left( \int_{\R^2}(|\nabla^2 v|^2 +
|\nabla^3 u|^2)dx \right)^{\frac{1}{2}}\\
& \ \cdot \left( \int_{\R^2} (|v|^4
+ |\nabla u|^4)dx \right)^{\frac{1}{2}}\\
\leq & \ \frac{a}{8}\int_{\R^2} (|\nabla^2 v|^2 + |\nabla^3
u|^2)dx + \left( C\int_{\R^4} (|v|^4 + |\nabla u|^4)dx
\right)\\
& \ \cdot \left( \int_{\R^2} (|\nabla v|^2 + |\nabla^2 u|^2)dx
\right).
\endsplit
$$

This, together with (3.25), shows that for $t\in(0,T)$,
$$\split & \frac{d}{dt} \int_{\R^2} (|\nabla v(x,t)|^2
+ |\nabla^2 u(x,t)|^2)dx + \frac{a}{8}\int_{\R^2} (|\nabla^2
v(x,t)|^2 + |\nabla^3 u(x,t)|^2)dx\\
\leq & \left (C\int_{\R^2} (|\nabla u|^4 + |v|^4)dx \right)
\int_{\R^2} (|\nabla v(x,t)|^2 + |\nabla^2 u(x,t)|^2) dx.
\endsplit\tag 3.26
$$

It follows from (3.9), (3.26), and Gronwall's inequality that for
any $s$ and $t$ with $\tau\leq s<t\leq T$,
$$\split
& \ \int_{\R^2} (|\nabla v|^2 + |\nabla^2 u|)
(x,t)\,dx\\
\leq & \ \left( e^{C\int^t_s\!\!\int_{\R^2} (|\nabla u|^4 + |v|^4)
(x,t)\,dx\,dl} \right) \left( \int_{\R^2} (|\nabla v|^2 +
|\nabla^2 u|^2)(x,s)\,dx\right)\\
\leq & \ \left( e^{C\,\varepsilon_1(1+TR^{-2})} \int_{\R^2}
e(u_0,v_0) dx \right) \int_{\R^2} (|\nabla v|^2 + |\nabla^2
u|^2)(x,s)dx
\endsplit\tag 3.27
$$
Thanks to (3.8), (3.27), and the mean value theorem, we conclude
that
$$\split
& \ \sup_{\tau\leq t\leq T} \int_{\R^2} (|\nabla^2 u|
+ |\nabla v|^2)(\cdot,t)dx\\
\leq & C\tau^{-1} (1+TR^{-2}) \, E(u_0,v_0) \,
e^{C\,\varepsilon_1(1+TR^{-2}) E(u_0,v_0)}
\endsplit\tag 3.28
$$
for any $\tau>0$. Then, by a similar proof as in Lemma 4, we can
show that $u$ belongs to $C^{1/8}(\R^2\times [\tau ,T]$ for any
$\tau >0$. In the appendix below (Section 5), we can show that
$(u,v)$ is regular for all $t\in (0,T]$. \qed
\enddemo

\proclaim {Remark} Let $u\in V(0 ,T)$ and $v\in H(0,T)$ be a
solution of (1.5)-(1.7) with initial values $u_0\in H_b^2(\R^2;
S^2)$, $v_0\in H^1(\R^2; \R^2)$ and $\text {div }v_0 =0$. Assume
that there are constants $\varepsilon_1$ and $R_0>0$ such that
$$\esssup_{0\leq t\leq T, x\in \R^2} \int_{B_R(x)}|\nabla
u(\cdot ,t)|^2+|v\left(\cdot , t\right)|^2 \,dx<\varepsilon_1
 $$ for  any
$R\in (0, R_0]$. Then, for $t\in [0, T]$, we have
$$\split
& \ \sup_{0\leq t\leq T} \int_{\R^2} (|\nabla^2 u(x,t) |^2+|\nabla
v (x,t)|^2)\,dx\\
\leq & \ C_3 (1 +T R^{-2})\, E(u_0,v_0)\,
e^{C\varepsilon_1(1+TR^{-2}) E(u_0,v_0)},
\endsplit\tag 3.29
$$
with $C_3=C(||u_0||_{H^2_b} + ||v_0||_{H^1})$.
\endproclaim

We are not able to prove the uniqueness of solutions to
(1.5)-(1.9) for initial value in $H^1\times L^2$ as one in \cite
{St1; Lemma 3.12}. However, we obtain

\proclaim {Lemma 11} Let $(u_1,v_1), (u_2,v_2)\in V(0 ,T)\times
H(0,T)$ be two smooth solutions of (1.5)-(1.7) with smooth initial
values $(u_0,v_0)\in H_b^2(\R^2; S^2)\times H^1(\R^2; \R^2)$ and
$\text {div }v_0 =0$. Then $(u_1,v_1)=(u_2,v_2)$.
\endproclaim
\demo {Proof} Following the proof of Proposition 15 in the
Appendix, we can assume that
$$ |\nabla u_1|+|\nabla u_2|+|v_1|+|v_2|\leq C
$$
for a constant $C>0$. For simplicity, we set in (1.7)
$$\split B(u,\nabla u):=&-W_{u^i}(u,\nabla
u)+W_{u^k}(u,\nabla u)u^k u^i+W_{p_{\a}^k}(u,\nabla
u)\nabla_{\a}u^k u^i \\
&+V_{p_{\a}^k}(u,\nabla u) u^k\nabla_{\a}u^i
\endsplit
$$
It follows from (1.7) that
$$\split
&\frac 1 2\frac d{dt}\int_{\R^2} | (u_1-u_2)|^2\,dx + \int_{\R^2}
(W_{p_{\a}^i}(u_1,\nabla u_1)-W_{p_{\a}^i}(u_2,\nabla u_2))\nabla_\a (u^i_1-u^i_2)\,dx\\
&=\int_{\R^2} (u_1^ju_1^iV_{p_{\a}^j}(u_1,\nabla u_1)-
u_2^ju_2^iV_{p_{\a}^j}(u_2,\nabla u_2))\nabla_\a
(u^i_1-u^i_2)\,dx\\
&\quad -\int_{\R^2}[-B(u_1,\nabla u_1) + B(u_2,\nabla u_2) +
(v_1\cdot\nabla u_1)-(v_2\cdot \nabla )u_2 ]\cdot (u_1-u_2)\,dx\\
& :=I_5+I_6.
\endsplit \tag 3.30
$$
By  Young's inequality, the last term on the right hand of the
above identity can be estimated as
$$\split &I_6=-\int_{\R^2}[-B(u_1,\nabla u_1) + B(u_2,\nabla u_2) +
(v_1\cdot\nabla u_1)-(v_2\cdot \nabla )u_2 ]\cdot
(u_1-u_2)\,dx\\
&\leq C\int_{\R^2} |u_1-u_2|^2+|v_1-v_2|^2\,dx +\frac a4
\int_{\R^2} |\nabla (u_1-u_2)|^2\,dx.
\endsplit
$$
The difficult part is to estimate $I_5$. Using an uniform open
ball covering of $\R^2$, we can estimate only the local integral
$$\int_{B_{r_0}(x_0)}(u_1^ju_1^iV_{p_{\a}^j}(u_1,\nabla u_1)-
u_2^ju_2^iV_{p_{\a}^j}(u_2,\nabla u_2))\nabla (u_1-u_2)\,dx.
$$

 Now we can think about in the equation (1.7) with $\frac {\partial u}{\partial x_3}=0 $ in a domain of $\R^3.$
After a rotation $\Cal R\in O(3)$, the integrand (1.2) has the
following invariant property:
$$W(\Cal R u, \Cal R \D u\Cal R^T)=W(u,\D u).$$
Therefore,  the system (1.5)-(1.7) is invariant for a rotation.
Without loss of generality, we can assume that $u_0(x_0)=(0,0,1)$.
Since $u_1$ and $u_2$ are uniformly continuous in $(x,t)\in
\R^2\times [0,\tau ]$ for some $\tau >0$, there exists a constant
$r_0
>0$ such that for any $(x,t)\in B_{r_0}(x_0)\times [0,\tau ]$
$$ |u_1(x,t)-u_0(x_0)|\leq \varepsilon ,\quad |u_2(x,t)-u_0(x_0)|\leq \varepsilon .
$$
Then
$$\split
&\int_{B_{r_0}(x_0)}(u_1^ju_1^iV_{p_{\a}^j}(u_1,\nabla u_1-\nabla
u_2)\nabla_\a (u^i_1-u^i_2)\,dx\\
&\leq C\varepsilon \int_{B_{r_0}(x_0)}|\nabla (u_1-u_2)|^2\,dx
+C\int_{B_{r_0}(x_0)}|\nabla (u^3_1-u^3_2)|^2\,dx.
\endsplit
$$
It follows from $|u|=1$ that
$$u^3\nabla_l u^3=-u^1\nabla_l u^1-u^2\nabla_l u^2.
$$
 Then an elementary calculation shows that
$$|\nabla (u^3_1-u^3_2)|\leq C |u_1-u_2| + C\varepsilon  |\nabla
(u_1-u_2)|.
$$
By a covering argument, we apply all above estimates to obtain
$$I_5\leq C\varepsilon \int_{\R^2}|\nabla (u_1-u_2)|^2\,dx
+C\int_{\R^2}|(u_1-u_2)|^2\,dx.
$$
Therefore, choosing $\varepsilon$ sufficiently small yields
$$\frac 1 2\frac d{dt}\int_{\R^2} | (u_1-u_2)|^2\,dx + \frac a 2\int_{\R^2}
|\nabla (u_1-u_2)|^2\,dx\leq
C\int_{\R^2}|(u_1-u_2)|^2+|v_1-v_2|^2\,dx\tag 3.31
$$

Using (1.5) and (1.6), one can obtain
$$\split
&\quad \frac 1 2 \frac d{dt}\int_{\R^2} |v_1-v_2|^2\,dx +
\int_{\R^2}
|\nabla (v_1-v_2)|^2\,dx\\
&\leq \tilde C\int_{\R^2} ( |v_1-v_2|^2 +|u_1-u_2|^2+|\nabla
(u_1-u_2)|^2) \,dx + \frac 1 2  \int_{\R^2} |\nabla
(v_1-v_2)|^2\,dx.
\endsplit
\tag 3.32
$$
Combining (3.31) with (3.32) gives
$$\split
&\quad \frac 1 2 \frac d{dt}\int_{\R^2} (\tilde C |u_1-u_2|^2+
\frac a 4 |v_1-v_2|^2)\,dx\leq C\int_{\R^2}  (\tilde C
|u_1-u_2|^2+ \frac a 4 |v_1-v_2|^2) \,dx.
\endsplit
\tag 3.33
$$
Integrating (3.33) in $t$ and applying the Gronwall inequality, we
conclude
$$\int_{\R^2} ( \tilde C |u_1-u_2|^2+
\frac a 4 |v_1-v_2|^2) (\cdot , t)\,dt \leq C\int_{\R^2} ( \tilde
C |u_1-u_2|^2+ \frac a 4 |v_1-v_2|^2)(\cdot , 0)\,dt=0.
$$
 This proves our claim. \qed
\enddemo

\head {\bf 4. Local existence and Proof of Theorem B}\endhead

In this section, we prove local existence of solutions of
(1.5)-(1.7) and complete the proof of Theorem B. Recall the
notation that $\tilde{V}(\tau,t)$ denotes the space $V(\tau,t)$
where $S^2$ is replaced by $\R^3$.

\proclaim {Lemma 12}   For a pair $(u_0, v_0)\in H_b^{1}(\R^2,
S^2)\times L^{2}(\R^2,\R^2)$ with $\text {div }v_0 =0$ in $\R^2$
in the sense of distribution,  there is a local regular solution
$(u_{\varepsilon}, v_{\varepsilon})\in \tilde{V}(0,T)\times
H(0,T)$ of (1.10)-(1.12) with initial data (1.9) for some $T>0$.
\endproclaim
\demo {Proof} Although Lin-Liu proved only the global existence of
the solution to (1.10)-(1.12) with initial data (1.9) for the case
of $k_1=k_2=k_3$, their proofs still work for the local existence
for the system (1.10)-(1.12). Thus we omit the details and refer
readers to \cite {LL1} and \cite {LL2}.\qed
\enddemo

\proclaim {Theorem 13} {\bf(Local existence)} For a pair $(u_0,
v_0)\in H_b^{1}(\R^2, S^2)\times L^{2}(\R^2,\R^3)$ with $\text
{div }v_0 =0$ in $\R^2$ in the sense of distribution, there is a
local solution $(u, v)\in V(0, t_1 )\times H(0,t_1)$ of
(1.5)-(1.7) with initial value $(u_0, v_0)$ for some $t_1>0$.
\endproclaim
\demo {Proof} For any map $u_0\in H_b^{1}(\R^2, S^2)$, one can
approximate it by a sequence of smooth maps in $H^{1}_b(\R^2,
S^2)$. Without loss of generality, we assume that $u_0\in
H_b^{2}(\R^2, S^2)$ and $v_0\in H^1(\R^2, \R^3)$ with $\text {div
}v_0 =0$ in $\R^2$  are smooth. Then thanks to Lemma 12, there is
a local regular solution $(u_{\varepsilon}, v_{\varepsilon})\in
\tilde{V}\times H$ of (1.10)-(1.12) with initial data (1.9).

For each pair $(u,v)$, set
$$e_{\varepsilon} (u,v) = W(u,\nabla u) + \frac 1 {2\varepsilon^2} (1-|u|^2)^2 + |v|^2,
\qquad E(u,v) = \int_{\R^2} e_\varepsilon (u,v) dx.$$

Then same calculations as for (3.1) give
$$ E(u_\varepsilon(\cdot,t) ,v_\varepsilon(\cdot,t)) +
2\int^t_0\!\!\int_{\R^2} (|\partial_t\,u_\varepsilon +
(v_\varepsilon\cdot\nabla u_\varepsilon)u_\varepsilon|^2 + |\nabla
v_\varepsilon|^2)^2\,dx\,dt = E(u_0,v_0).\tag 4.1
$$

By a similar analysis as in the proof of Lemma 8 and Lemma 7, one
can show that there exist uniform positive constants $R_0$ and
$\varepsilon_1$, and a positive time
$T_{\varepsilon}=T(\varepsilon, R_0, \varepsilon_1)$ such that the
problem (1.10)-(1.12) with initial data (1.9) has a regular
solution $(u_\varepsilon,
v_\varepsilon)\in\tilde{V}(0,T_\varepsilon)\times
H(0,T_\varepsilon)$ for each fixed $\varepsilon>0$, and
furthermore, it holds that
$$\sup_{0\leq t\leq T_{\varepsilon}} \int_{B_R(x_0)}|\nabla u_{\varepsilon}(\cdot ,t) |^2
+|v_{\varepsilon}(\cdot ,t) |^2
+\frac 1{2\varepsilon^2}  (1-|u_{\varepsilon}(\cdot ,t) |^2)^2
\,dx<\varepsilon_1\tag 4.2
$$
for any positive $R\leq R_0$.

Next, we will show that there is a constant $t_1>0$, independently
of $\varepsilon$, such that $T_{\varepsilon}\geq t_1$ and the
solutions $(u_{\varepsilon}, v_{\varepsilon})$ is bounded in
$\tilde{V}(0, t_1) \times H (0, t_1)$ uniformly in $\varepsilon$.

First, we claim that for all $t\in [0,\,\min\{1,
T_{\varepsilon}\}]$
$$\frac 12 \leq |u_{\varepsilon}(x,t)|\leq \frac 32.\tag 4.3
$$

To verify (4.3), we re-scale the solution by
$$\tilde
u(x,t)=u_{\varepsilon} (\varepsilon x, \varepsilon^2 t),\quad
\tilde v(x,t)={\varepsilon } v_{\varepsilon} (\varepsilon x,
\varepsilon^2 t), \quad \tilde P(x,t)= {\varepsilon^2 }
P_{\varepsilon} (\varepsilon x, \varepsilon^2 t).$$ Then $(\tilde
u,\tilde v)$ solves the following approximate Ericksen-Leslie
system

$$\tilde v^i_t+(\tilde v\cdot\nabla )\tilde v^i- \laplacian \tilde v^i+\nabla_{x_i} \tilde P=
-  \nabla_{x_j}(\nabla_{x_i}\tilde u^k\,W_{p_j^k}(\tilde u,\nabla
\tilde u)), \tag 4.4
$$
$$ \nabla \cdot \tilde v=0,\tag 4.5
$$
$$
 \tilde u^i_t+(\tilde v\cdot \nabla  )\tilde u^i =   \nabla_{\a} \left [W_{p_{\a}^i} (\tilde u,
\nabla \tilde u) \right ]-W_{u^i}(\tilde u,\nabla \tilde u)+
\tilde u^i (1-|\tilde u |^2) \tag 4.6
$$
for $i=1,2,3$, with initial data
$$\tilde v(x,0)={\tilde v}_0(x), \quad \tilde u(x,0)={\tilde u}_0(x), \quad \forall x\in \R^2 , \tag
4.7
$$
where $\tilde u_0(x)=u_0(\varepsilon x)$ and  $\tilde v_0 (x)=
 \varepsilon v_0(\varepsilon x)$ satisfy
$$\int_{\R^2} e( \tilde
u_0(x), \tilde v_0(x))\,dx =\int_{\R^2} e(u_0(x), v_0(x))\,dx.$$
The condition (4.2) becomes
$$\esssup_{0\leq t\leq  \frac {T_{\varepsilon} }{\varepsilon^2}, x\in \R^2}
 \int_{B_{\frac R{\varepsilon }}(x)}|\nabla
\tilde u(\cdot ,t)|^2+\frac 1 2 (|1-|\tilde u (\cdot
,t)|^2)^2+|\tilde v|^2\,dx<\varepsilon_1\tag 4.8
$$
for any $R\in (0, R_0]$. While the basic energy identity (4.1)
becomes
$$\split
& \ \int_{\R^2} \left[
W(\tilde{u},\nabla\tilde{u})+\frac{1}{2}(1-|\tilde{u}|^2)^2
\right](\cdot,t)dx\\
& \ + \int^t_0\!\!\int_{\R^2} (|\partial_t\,\tilde{u}|^2 +
|\partial_t\,\tilde{u} +
(\tilde{v}\cdot\nabla)\tilde{u}|^2)(\cdot,t) dx\,dt = E(u_0,v_0)
\endsplit\tag 4.9
$$
for all $t\in(0,\frac{T_\varepsilon}{\varepsilon^2})$.

Without loss of generality, we assume $T_{\varepsilon}\leq 1$. Let
$\tau$ be the maximal time in $[0,\frac
{T_{\varepsilon}}{\varepsilon^2} ]$ such that
$$\frac 1 4 \leq |\tilde u (x,t)|\leq 2.\tag 4.10
$$
By (4.1) and similar arguments as for Lemma 8, one can derive from
(4.4)-(4.6) that there exists a uniform constant $C_0$ such that

$$\split
&\quad C_0 \int_0^{\tau}\int_{\R^2} \left (
 |\nabla^2\tilde u |^2+|\nabla\tilde v|^2\right ) \,dx\,dt\\
&\leq -\int_0^{\tau }\int_{\R^2}\laplacian \tilde u\cdot \tilde u
(1-|\tilde u |^2)\,dx\,dt + C \int_{\R^2}|\nabla
u_0|^2 +|v_0|^2 \, dx\\
& + C\int_0^{\tau }\int_{\R^2} |\nabla \tilde u|^4 +|\tilde
v|^4\,dx\,dt .
\endsplit\tag 4.11
$$

Integration by parts yields
$$\split
& \quad -\int_0^{ \tau}\int_{\R^2}\laplacian \tilde u \
 \cdot\tilde  u (1-|\tilde u |^2)\,dx\,dt\\
 &=\int_0^{\tau}\int_{\R^2}|\nabla
 \tilde u |^2 (1-|\tilde u |^2)\,dx\,dt -\frac 1
 2\int_0^{\tau}\int_{\R^2}|\nabla |\tilde u
 |^2|^2\,dxdt
 \endsplit\tag 4.12
$$
By Young's inequality and using (4.1), (4.6) and (4.10), one can
obtain
$$\split
&\int_0^{\tau}\int_{\R^2}|\nabla
 \tilde u |^2 (1-|\tilde u |^2)\,dx\,dt\leq \eta \int_0^{\tau}\int_{\R^2} (1-|\tilde u
|^2)^2\,dx\,dt + C\int_0^{\tau}\int_{\R^2}|\nabla
 \tilde u |^4\,dx\,dt\\
 & \leq \eta \int_0^{\tau}\int_{\R^2} |\nabla^2 u|^2\,dx\,dt
+ C\int_0^{\tau}\int_{\R^2}|\nabla
 \tilde u |^4\,dx\,dt +C\int_{\R^2}(|\nabla u_0|^2 +|v_0|^2)\,dx
 \endsplit\tag 4.13
$$
for a small constant $\eta$.

Combining (4.11)-(4.13) and choosing $\varepsilon_1$ sufficiently
small in (4.2) with Lemma 2, we conclude that
$$ \int_0^{\tau}\int_{\R^2}
\left ( |\nabla^2\tilde u |^2+|\nabla\tilde v|^2\right ) \,dx\,dt
+ \int^\tau_0\!\!\int_{\R^2} (1-|\tilde{u}|^2)^2\,dx\,dt \leq C(1+
\frac { \tau \varepsilon^2 } {R^2}) \int_{\R^2}e( u_0,v_0)\,dx\tag
4.14$$ for any $R\leq R_0$.

It follows also from Lemma 2, (4.8), (4.9), and (4.14) that
$$\int^\tau_0\!\!\int_{\R^2} (|\nabla\tilde{u}|^4 + |\tilde{v}|^4)\,dx\,dt \leq
C\,\varepsilon_1 \left( 1 + \frac{\tau\varepsilon^2}{R^2}\right)
\, E(u_0,v_0)\tag 4.15
$$
for any $R\leq R_0$.

Now following the calculation for (3.25), one can derive that for
any $t\in(0,\tau)$,
$$\split
& \ \frac d{dt} \left (\int_{\R^2} (|\nabla^2\,\tilde{u}|^2 +
|\nabla \tilde{v}|)(\cdot,t)dx \right ) + \frac{a}{4} \int_{\R^2}
(|\nabla^3 \tilde{u}|^2
+ |\nabla^2 \tilde{v}|)(\cdot,t)dx \\
\leq & \ C\int_{\R^2} (|\tilde v|^2 +|\nabla \tilde u|^2) (|\nabla
\tilde v|^2 + |\nabla^2 \tilde u|^2)(\cdot,t)dx\\
& \ + C\left| \int_{\R^2} \nabla_\beta (\tilde{u}^i
(1-|\tilde{u}|^2)) \cdot \nabla_\beta\,\laplacian\tilde{u}^i\,dx
\right| + C\int_{\R^2} |\nabla\tilde{u}|^6\,(\cdot,t)dx.
\endsplit\tag 4.16
$$

Note that
$$\split
& \ C\left| \int_{\R^2} \nabla_\beta (\tilde{u}^i
(1-|\tilde{u}|^2)) \cdot
\nabla_\beta\,\laplacian\tilde{u}^i\,dx\right|\\
= & \ C\left| \int_{\R^2} \laplacian (\tilde{u}^i
(1-|\tilde{u}|^2)) \cdot
\laplacian\tilde{u}^i\,dx \right|\\
\leq & \ C\int_{\R^2} |\nabla^2\tilde{u}(\cdot,t)|^2 dx +
C\int_{\R^2} |\nabla\tilde{u}(\cdot,t)|^4\,dx
\endsplit\tag 4.17
$$
and
$$\split
C\int_{\R^2} |\nabla\tilde{u}(\cdot,t)|^6\,dx = & \ -C\int_{\R^2}
|\nabla\tilde{u}(\cdot,t)|^4\,\tilde{u} \cdot
\laplacian\tilde{u}\,dx - C\int\nabla_\a (|\nabla\tilde{u}|^4)
\tilde{u}\cdot\nabla_\a\,\tilde{u}\,dx\\
\leq & \ \frac{C}{2}\int_{\R^2} |\nabla\tilde{u}(\cdot,t)|^6\,dx +
C\int_{\R^2} |\nabla\tilde{u}(\cdot,t)|^2 \,
|\nabla^2\tilde{u}(\cdot,t)|dx.
\endsplit\tag 4.18
$$

It follows from (4.16)-(4.18) that
$$\split
& \frac{d}{dt} \left( \int_{\R^2} (|\nabla^2\tilde{u}|^2 +
|\nabla\tilde{v}|^2)(\cdot,t)dx \right) + \frac{a}{4} \int_{\R^2}
(|\nabla^3\tilde{u}|^2 + |\nabla^2\tilde{v}|^2)(\cdot,t)dx\\
\leq & C\int_{\R^2} (|\tilde{v}|^2 + |\nabla\tilde{u}|^2)
(|\nabla\tilde{v}|^2 + |\nabla^2\tilde{u}|^2)dx + C \left(
\int_{\R^2} (|\nabla^2\tilde{u}(\cdot,t)|^2 +
|\nabla\tilde{u}(\cdot,t)|^4)dx\right)
\endsplit
$$

Using the Gagliardo-Nirenberg-Sobolev's inequality, one can get
$$\split
& \ C\int_{\R^2} (|\tilde{v}|^2 + |\nabla\tilde{u}|^2)
(|\nabla\tilde{v}|^2 + |\nabla^2\tilde{u}|^2)dx\\
\leq & \ C \left( \int_{\R^2} (|\tilde{v}|^4 +
|\nabla\tilde{u}|^4)dx \right)^{\frac{1}{2}} \left( \int_{\R^2}
(|\nabla\tilde{v}|^4 + |\nabla^2\tilde{u}|^4)dx
\right)^{\frac{1}{2}}\\
\leq & \ C\left( \int_{\R^2} (|\nabla\tilde{v}|^2 +
|\nabla^2\tilde{u}|^2)dx \right)^{\frac{1}{2}} \left( \int_{\R^2}
(|\nabla^2\tilde{v}|^2 + |\nabla^3\tilde{u}|^2)dx
\right)^{\frac{1}{2}}\\
& \ \cdot \left( \int_{\R^2} (|\tilde{v}|^4 +
|\nabla\tilde{u}|^4)dx \right)^{\frac{1}{2}}\\
\leq & \ \frac{a}{8}\int_{\R^2} (|\nabla^2\tilde{v}|^2 +
|\nabla^3\tilde{u}|^2)dx + \left( C\int_{\R^2} (|\tilde{v}|^4 +
|\nabla\tilde{u}|^4)dx \right)\\
& \ \cdot \int_{\R^2} (|\nabla\tilde{v}|^2 +
|\nabla^2\tilde{u}|^2)dx.
\endsplit
$$

This, together with (4.18), shows that
$$\split
& \ \frac{d}{dt}\int_{\R^2} (|\nabla^2\tilde{u}|^2 +
|\nabla\tilde{v}|^2)(\cdot,t)dx + \frac{a}{8}\int_{\R^2}
(|\nabla^3\tilde{u}|^2 + |\nabla^2\tilde{v}|)(\cdot,t)dx\\
\leq & \ \left( C\int_{\R^2} (|\tilde{v}|^4 +
|\nabla\tilde{u}|^4)(\cdot,t)dx \right) \int_{\R^2}
(|\nabla^2\tilde{u}|^2 + |\nabla\tilde{v}|^2)(\cdot,t)dx + h(t)
\endsplit\tag 4.19
$$
with
$$
h(t)=C\int_{\R^2} (|\nabla^2\tilde{u}|^2 +
|\nabla\tilde{u}|^4)(\cdot,t)dx.\tag 4.20
$$

It then follows from (4.14), (4.15), (4.19)-(4.20), and Gronwall's
inequality that for all $t\in(0,\tau)$,
$$\split
& \ \int_{\R^2} (|\nabla^2\tilde u|^2 +
|\nabla\tilde{v}|^2)(\cdot,t)dx\\
\leq & \ e^{C\int^t_0\!\!\int_{\R^2} (|\tilde{v}|^4 +
|\nabla\tilde{u}|^4)(\cdot,l)dx\,dl} \cdot \int_{\R^2} (|\nabla^2
u_0|^2 + |\nabla v_0|^2)dx \\
& \ + \int^t_0 e^{C\int^t_s\!\!\int_{\R^2} (|\nabla\tilde{u}|^4 +
|\tilde{v}|^4)(\cdot,l)dx\,dl} \cdot h(s)\,ds\\
\leq & \
e^{C\varepsilon_1(1+\frac{\tau\varepsilon^2}{R^2})\,E(u_0,v_0)}
\left( ||u_0||^2_{H^2} + ||v_0||^2_{H^1} + \int^t_0 h(s)ds
\right)\\
\leq & \ \left( ||u_0||^2_{H^2} + ||v_0||^2_{H^1} + C\left(
1+\frac{\tau\varepsilon^2}{R^2}\right) \, E(u_0,v_0) \right) \,
e^{C\varepsilon_1(1+\frac{\tau\varepsilon^2}{R^2})\,E(u_0,v_0)}.
\endsplit\tag 4.21
$$

Suppose that there is a $x_1\in \R^2$ such that $|\tilde
u(x_1,t)|<1/2$ (or $|\tilde u(x_1,t)|> \frac 3 2 $) with some
$t\in [0, \tau ]$. It follows from (4.21) that $\tilde{u}$ is
$C^{\frac 1 8}$-continuous uniformly in $(x,t)$. Then, there is a
constant $C_4$ so that for $x\in B_{1/4C_4}(x_1)$ with $\frac
{1}{4C_4}< \frac {R_0}{\varepsilon}$, we have
$$(1-|\tilde u(x,t)|^2)^2 \geq \frac 12.
$$
Then
$$
\frac 1 2 \int_{B_{\frac {1}{4C_4}}(x_1)}(1-|\tilde u
(x,t)|^2)^2\,dx \geq \frac 1 8|B_{\frac 1 {4C_4}}(0)|
>\varepsilon_1
$$
which contradicts (4.8) for a sufficiently small $\varepsilon_1$.
This shows that   our claim (4.3) holds for all $t\in [0,\frac
{T_{\varepsilon}}{\varepsilon}]$.

Finally, we show that $(u_{\varepsilon},v_{\varepsilon})$ is
bounded in $\tilde{V} (0, \min \{T_{\varepsilon} ,1\})\times H(0,
\min \{T_{\varepsilon} ,1\})$ uniformly for any positive
$\varepsilon<4C_4R_0$.

For any $t\leq\min(1,T_\varepsilon)$, it follows from (4.14) and
(4.15) that
$$\split
& \ \int^t_0\!\!\int_{\R^2} (|\nabla^2\,u_\varepsilon|^2 + |\nabla
v^\varepsilon|^2)(x,t) dx\,dt + \frac{1}{4\varepsilon^4}
\int^t_0\!\!\int_{\R^2} (1-|u_\varepsilon|^2)^2\,dx\,dt\\
\leq & \ C\left( 1+\frac{t}{R^2}\right)\,E(u_0,v_0),
\endsplit\tag 4.22
$$
$$ \int^t_0\!\!\int_{\R^2} (|\nabla u_\varepsilon|^4 +
|v_\varepsilon|^4)(x,t) dx\,dt \leq C\,\varepsilon_1 \left(
1+\frac{t}{R^2}\right)\,E(u_0,v_0).\tag 4.23
$$

Let $\varphi$ be the cut-off function as in the proof of Lemma 9.
Then by a similar analysis as in (3.14)-(3.17) and using
(4.22)-(4.23), one can get
$$\split
& \ \int_{\R^2} e(u_\varepsilon (\cdot,t), v_\varepsilon(\cdot,t))
\varphi^2\,dx + \int^t_0\!\!\int_{\R^2} (|\nabla v_\varepsilon|^2
+ \frac{1}{2}|\partial_t\,u_\varepsilon + v_\varepsilon\cdot\nabla
u_\varepsilon|^2) \varphi^2\,dx\,dt\\
\leq & \ \int_{\R^2} e(u_0,v_0) \varphi^2\,dx +
C\frac{t^{\frac{1}{2}}}{R} \left( 1+\frac{t}{R^2}
\right)^{\frac{1}{2}} \, E(u_0,v_0)\\
& \ + \int^t_0\!\!\int_{\R^2} \frac{1}{2\varepsilon^2}
(1-|u_\varepsilon|^2)^2\,
|v_\varepsilon\cdot\nabla(\phi^2)|dx\,dt.
\endsplit\tag 4.24
$$

On the other hand,
$$\split
& \ \int^t_0\!\!\int_{\R^2} \frac{1}{2\varepsilon^2}
(1-|u_\varepsilon|^2)^2 \,
|v_\varepsilon\cdot\nabla(\varphi^2)|\,dx\,dt\\
\leq & \ C\left( \int^t_0\!\!\int_{\R^2} \frac{1}{4\varepsilon^4}
(1-|u_\varepsilon|^2)^4\,dx\,dt \right)^{\frac{1}{2}} \left(
\int^t_0\!\!\int_{\R^2} \frac{|v_\varepsilon|^2}{R^2} dx\,dt
\right)^{\frac{1}{2}}\\
\leq & \ C \frac{t^{\frac{1}{2}}}{R} \left( \int_{\R^2}
\frac{1}{4\varepsilon^4} (1-|u_\varepsilon|^2)^2\,dx\,dt
\right)^{\frac{1}{2}} (E(u_0,v_0))^{\frac{1}{2}}\\
\leq & \ C \frac{t^{\frac{1}{2}}}{R} \left( 1+\frac{t}{R^2}
\right)^{\frac{1}{2}}\, E(u_0,v_0)
\endsplit\tag 4.25
$$
where one has used (4.1) and (4.22).

Hence,
$$ \int_{B_R(x_0)} e_\varepsilon (u_\varepsilon(x,t),
v_\varepsilon(x,t))\,dx \leq \int_{B_{RR}(x_0)} e(u_0,v_0)dx +
C\frac{t^{\frac{1}{2}}}{R}
\left(1+\frac{t}{R^2}\right)^{\frac{1}{2}} \, E(u_0,v_0) \tag 4.26
$$
for any $R\leq R_0$. First, choosing $R_1>0$ so that
$$\int_{B_{2R_1}(x_0)} e(u_0,v_0)dx < \frac{\varepsilon_1}{2}\tag
4.27
$$
for all $x_0\in\R^2$. Then, set
$$ t_1=\min\left\{ R_1,
\frac{\varepsilon_1\,R_1}{4C(\varepsilon_1+ E(u_0,v_0))}.
\right\}\tag 4.28
$$

Then for $t\leq t_1$,
$$\int_{B_R(x_0)} e_\varepsilon (u_\varepsilon(x,t),
v_\varepsilon(x,t))dx <\varepsilon_1 $$ for all $x_0\in\R^2$ and
$R\leq R_1$. Consequently, we have shown that there is a uniform
$t_1\leq\min\{T_\varepsilon, 1\}$ such that $(u_\varepsilon,
v_\varepsilon)$ is bounded in $\tilde{V}(0,t_1)\times H(0,t_1)$
with $t_1$ independent of $\varepsilon$. Letting
$\varepsilon\rightarrow 0$, we can prove the local existence of
solution $(u,v)\in V(0,t_1)\times H(0,t_1)$ with initial data
(1.9). \qed
\enddemo

Now we complete the proof of Theorem B.

\demo{Proof of Theorem B} By Theorem 13, there is a local solution
$(u, v)\in V(0,t_1)\times H(0,t_1)$ of (1.5)-(1.7) in $\R^2\times
[0, t_1]$ with initial conditions (1.9) for some $t_1>0$. By
Lemmas 10-11, the solution can be extended in $[0,T_1)$ for a
maximal times $T_1$ such that at $T_1$, there is at least a
singular point $ x_{i}^{1} \in\R^2 $
 such that
$$\limsup_{t \nearrow  T_{1}} \int_{B_{R} ( x_{i}^{1})}e (
u,v)(\cdot , t) \,dx \geq \varepsilon_{0}$$
 for any  $R\leq R_0$ for some $R_0>0$ and $\varepsilon_0>0$.
 It is easy to see the solution $(u, v)\in V\times H$ is  regular for all $t\in (0, T_1)$.
 Then there exists a sequence of $\{t_n\}$ such that the sequence
 $(u(t_n),v(t_n))$ converges weakly to $(u(T_1), v(T_1))$ in
 $H^1(\R^2; S^2)\times L^2(\R^2; \R^3)$
 satisfying
 $$\int_{\R^2} e(u(T_1), v(T_1))\,dx\leq \int_{\R^2} e(u_0,
 v_0)\,dx-\varepsilon_0,\quad \text{div }v(T_1)=0.
 $$
Using the energy identity,    there is a finite
 number of  singular times $\{T_l\}_{l=1}^L$ in Theorem B.
\qed
\enddemo

\head {\bf 5. Appendix: The liquid crystal flow and regularity
issue}\endhead
 In this section, we  formulate the liquid crystal heat flow and discuss $C^{1,\a}$-regularity issues
 for  solutions of the liquid crystal flow (1.4) and the system (1.5)-(1.7).

    The liquid
 crystal equilibrium system in a form of vectors and tensors was derived by Hardt, Lin and
 Kinderlehler in \cite {HLK} using the Lagrange multiplier method,
 but we need a precise form of (1.3) in coordinates.

Let $\phi$ be a smooth functional in $C_0^{\infty}(\Omega ,
\R^3)$. We consider a variation
$$u_t(x)=\frac {u+t\phi }{|u(x)+t\phi(x) |}=\frac {u+t\phi }{(1+2tu\cdot \phi +t^2\phi^2)^{1/2}}
$$
and compute
$$\frac {du_t}{dt}=\phi -\frac {(u+t\phi )(u\cdot \phi +t|\phi|^2)} {(1+2tu\cdot \phi
+t^2\phi^2)^{1/2}}.
$$
To derive the Euler-Lagrange equations, we compute
$$\left .\frac {d}{dt}\int_{\Omega} W(u_t,\nabla u_t)\,dx \right
|_{t=0}=0.
$$
This implies
$$\int_{\Omega} \left .\left (W_{u^j}\frac {du_t^j}{dt} +W_{p_{\a}^i}  \frac {d \nabla_{\a} u_t^i}{dt}\right )\right
|_{t=0}\,dx =0,
$$
where $W_{p_{\a}^i}(u,p)=\frac {\partial W}{\partial p_{\a}^i}$
and $W_{u^i}=\frac {\partial W}{\partial u^i}$. Note
$$\left .\frac {du_t^i}{dt}\right |_{t=0}=\phi^i-u^i  (u\cdot \phi
),\quad \left .\frac {d \nabla_{\a} u_t^i}{dt}\right
|_{t=0}=\nabla_{\a}\phi^i -\nabla_{\a}u^i (u\cdot \phi
)-u^i\nabla_{\a} (u\cdot \phi ).
$$
We conclude that
$$\split
\int_{\Omega}& W_{u^j}(u,\nabla u) \left [\phi^j-u^j (u\cdot \phi
)\right ]\\
&+W_{p_{\a}^i}(u,\nabla u) \left [ \nabla_{\a}\phi^i
-\nabla_{\a}u^i (u\cdot \phi )-u^i\nabla_{\a} (u\cdot \phi )\right
]  \,dx =0
\endsplit\tag 5.1
$$
for any $\phi\in  C_0^{\infty}(\Omega , \R^3)$. Therefore, we call
that $u\in H^1(\Omega , S^2)$ is a weak solution to the liquid
crystal system if $u$ satisfies
$$\split
&-\nabla_{\a} \left [W_{p_{\a}^i} (u,\nabla u)- u^ku^i
W_{p_{\a}^k}(u,\nabla u) \right ]+W_{u^i}(u,\nabla
u)\\
&-W_{u^k}(u,\nabla u)u^k u^i-W_{p_{\a}^k}(u,\nabla
u)\nabla_{\a}u^k u^i -W_{p_{\a}^k}(u,\nabla u) u^k\nabla_{\a}u^i=0
\endsplit
$$
in the sense of distribution. Note $|u|^2=1$, then $u^i\nabla
u^i=0$. This system  is the exact form of (1.3).

Then, the liquid crystal flow can be formulated as in (1.4), i.e.,
$$\split
\frac {\partial u^i}{\partial t}=&\nabla_{\a} \left [W_{p_{\a}^i}
(u,\nabla u)- u^ku^i V_{p_{\a}^k}(u,\nabla u) \right
]-W_{u^i}(u,\nabla
u)\\
&+W_{u^k}(u,\nabla u)u^k u^i+W_{p_{\a}^k}(u,\nabla
u)\nabla_{\a}u^k u^i +V_{p_{\a}^k}(u,\nabla u) u^k\nabla_{\a}u^i.
\endsplit
$$

Next, we will prove that a H\"older continuous solution of (1.4)
belongs to  $C^{1,\a}$ for some $\a$ with $0<\a <1$. For any point
$z_0=(x_0,t_0)\in \Omega \times [0,\R)$ and any number $R>0$, we
use standard notations:
$$\split
&B(x_0,R)=\{x\in\R^3: |x-x_0|<R\},\quad Q(z_0,R)=B(x_0,R)\times
(t_0-R^2,t_0),\\
& S_R(z_0)=B(x_0,R)\times \{t_0-R^2\}\cap \partial B(x_0,R)\times
(t_0-R^2, t_0).
\endsplit
$$

\proclaim {Proposition 14} Let $\Omega$ be a domain in $\R^3$ with
smooth boundary $\partial \Omega$. Let $u$ be  a weak solution of
(1.4) and H\"older continuous in $\Omega\times [0,T)$. Then,
$\nabla u$ is (locally) H\"older continuous with the same exponent
in $\Omega\times [0,T)$.
\endproclaim
\demo{Proof} Assume that $u(x,t)$ is H\"older continuous with
exponent $\b$, $0<\b<1$. Let $(x_0,t_0)\in \Omega\times (0,T)$
with $Q_{4R_0}(z_0)\subset \Omega\times (0,T)$ for some $R_0>0$.
Note $u(x_0,t_0)=e\in S^2$. After a rotation,  we can assume that
$e=(0,0,1)$.

It follows from $|u|=1$ and Cauchy's inequality that
$$|u^3|^2 \, |\D u^3|^2\leq (1-|u^3|^2) |\D u|^2\leq 2 |u-u(x_0,t_0) | |\D u|^2.\tag
5.2
$$

Denote
$$\tilde p =(p_{\a}^j)_{3\times 2}.
$$
 Using the structure  of $W(u, p)$, we can write
$$W_{\tilde p}(u, \D u)=\tilde W_{\tilde p} (u, \D  u^1, \nabla
u^2)+f(u,\D u^3),
$$
where $|f(u, \D  u^3)|\leq C|\D  u^3|$.

 Let $\tilde v= (v^1,  v^2 )$ be the solution of
  the Cauchy-Dirichlet problem
$$\split &v^i_t=\nabla_{\a} \left [W_{\tilde p_{\a}^i}
(e,\nabla v^1,\nabla v^2)\right ]\quad \text{in } Q_R(z_0)\\
&v^i=u^i\quad \text{on } S_R(z_0).
\endsplit\tag 5.3
$$
for $i=1,2$. Since (5.3) is a parabolic system with constant
coefficients, it follows from Proposition 1.2 in \cite{GS;
Proposition 1.2} that for all $\rho \leq R\leq R_0$
$$\int_{Q_{\rho}} |\nabla \tilde v|^2\,dz\leq C \left (\frac {\rho} R\right )^{5}\int_{Q_R} |\nabla
\tilde v|^2\,dz
$$
and
$$ \int_{Q_{\rho}} |\nabla \tilde v-(\nabla \tilde v)_{\rho} |^2\,dz\leq C \left (\frac {\rho} R\right )^{7}\int_{Q_R} |\nabla
\tilde v-(\nabla \tilde v)_{R}|^2\,dz.$$

 Set $\tilde w=\tilde u-\tilde v$. Then for all
$\rho <R$, we have
$$\int_{Q_{\rho}} |\nabla u|^2\,dz\leq C \left (\frac {\rho} R\right )^{5}\int_{Q_R} |\nabla
u|^2\,dz +C\int_{Q_R} |\nabla \tilde w|^2\,dz +C\int_{Q_R} |\nabla
u^3|^2\,dz \tag 5.4
$$
and
$$\split
\int_{Q_{\rho}} |\nabla u-(\nabla u)_{\rho}|^2\,dz\leq &C \left
(\frac {\rho} R\right )^{7}\int_{Q_R} |\nabla u-(\nabla
u)_R|^2\,dz +C\int_{Q_R} |\nabla \tilde w|^2\,dz \\
&+C\int_{Q_R} |\nabla u^3|^2\,dz. \endsplit\tag 5.5
$$
Note that $u$ is $\b$-H\"older continuous in $\Omega\times [0,T)$
and $u(x_0,t_0)= (0,0,1)$.

Although there is no maximum principle for the parabolic system
(5.3) with constant coefficients, Giaquinta-Struwe in \cite {GS;
page 445} obtained that
$$\sup_{Q_R} |v-u(x_0,t_0)|\leq C\sup_{Q_R} |u-u(x_0,t_0)|
$$with a constant $C$ independent of $R$ and $u$.
This implies
$$|\tilde w|\leq |u-u(x_0,t_0)|+|v-u(x_0,t_0)|\leq
CR^{\b}.\tag 5.6
$$
Multiplying the difference between (5.3) and (1.4) by $\tilde w^i$
($i=1,2$) and integrating over $Q_R$ lead to
$$\split
&  \int_{B_R} |\tilde w|^2(\cdot , t_0)\,dx+\int_{Q_R}
\sum_{i=1}^2 \nabla_{\a}\tilde w ^i W_{\tilde p_{\a}^i} (e,\nabla
\tilde
w)\,dx\\
&\leq \int_{Q_R} \sum_{i=1}^2 |\nabla_{\a}\tilde w^i|
\,|\tilde{W}_{\tilde p_{\a}^i} (e,\nabla \tilde
v)-\tilde{W}_{\tilde p_{\a}^i} (u,\nabla \tilde
u)|\,dx+C\int_{Q_R} |\nabla u^3||\nabla \tilde w|\,dx \\
&+\int_{Q_R} \sum_{i=1}^2\sum_{k=1}^3 \nabla_{\a}\tilde w^i u^i
u^k V_{p_{\a}^k} (u,\nabla u)+ C\int_{Q_R} |\tilde w| |\nabla
u|^2\,dx.
\endsplit \tag 5.7$$
Since $u$ is $\b$-H\"older continuous and $u(x_0,t_0)=(0,0,1)$, we
have $|u^i|\leq CR^{\b}$ for $i=1, 2$. Applying Young's inequality
and (5.2) yields
$$
\int_{Q_R} |\nabla \tilde w|^2\,dz\leq CR^{\b}\int_{Q_R}|\nabla
u|^2\,dz.\tag 5.8
$$
It follows that for all $\rho <R$,
$$\int_{Q_{\rho}} |\nabla u|^2\,dz\leq C \left (\frac {\rho} R\right )^{5}\int_{Q_R} |\nabla
u|^2\,dz +CR^{\b}\int_{Q_R}|\nabla u|^2\,dz.\tag 5.9
$$

We claim  the following Cacciopoli's inequality
$$\int_{Q(z_0,R)}|\nabla u|^2\,dz\leq C\frac 1 {R^2}\int_{Q(z_0,2R)} |u-  u_{2R}|^2\,dz\leq
CR^{3+2\b}.\tag 5.10
$$
for any $z_0\in\Omega\times (0,\infty )$ and $R\leq R_0$, where
$u_{2R}$ is the average of $u$ in $Q_{2R}(x_0,t_0)$.

Next, we prove this claim.  Let $\xi$ be a cut-off function in
$C_0^{\infty} (B_{2R}(x_0))$ with $0\leq \xi\leq 1$, $\xi\equiv 1$
in $B_{R}(x_0)$ and $|\nabla \xi|\leq \frac CR$. Let $\tau\in
C^{\infty}(\R, \R)$ be a function depends only on $t$ with $0\leq
\tau\leq 1$, $\tau\equiv 1 $ on $[t_0-R^2, t_0]$ and $\tau\equiv
0$ on $(-\infty , t_0-4R^2)$ and $|\partial_t \tau|\leq C/R^2$.

Testing (1.4) with $\phi =(u^i- u^i_{2R})\xi^2 \tau^2 I_{(-\infty
,t_0)}$ for $i=1,2$, where $I_{(-\infty ,t_0)}$ is the
characteristic function of  $(-\infty ,t_0)$, we have
$$\split
 &\int_{B_{2R}(x_0)}|u(\cdot , t_0)-{u}_{2R}|^2 \xi^2
\tau^2 (t_0)\,dx+  \int_{Q_{2R}(z_0)}
\sum_{i,j=1}^2W_{p_{\a}^ip_{\b}^j}
\nabla_{\a} u^i \nabla_{\a} u^j\xi^2 \tau^2\,dz\\
&\leq 2\int_{Q_{2R}(z_0)}\left [W_{p_{\a}^i} (u,\nabla u)- u^ku^i
V_{p_{\a}^k}(u,\nabla u) \right
]\nabla_{\a} \xi (u^i-  u_{2R}^i)\xi \tau^2 \,dz\\
&+C\int_{Q_{2R}(z_0)}|\nabla u|^2 |u- u_{2R} |\xi^2
\tau^2\,dz +2\int_{Q_{2R}(z_0)}|u-  u_{2R}|^2 \xi^2 \tau \partial_t \tau \,dz\\
& +\int_{Q_{2R}(z_0)}u^ku^i V_{p_{\a}^k} (u,\nabla u)
\,\nabla_{\a}( u^i- u^i_{2R}) \xi^2 \tau^2\,dz
+C\int_{Q_{2R}(z_0)}|\nabla u^3|^2\xi^2 \tau^2\,dz.
\endsplit
$$
Since $u$ is $\b$-H\"older continuous and $u(x_0,t_0)=(0,0,1)$,
$u(x,t)- u_{2R}$ can be chosen sufficiently small when $R_0$ is
small and $|u^1|+|u^2|$ is also small. We need to deal with the
above last term. By (5.2), the term $|\nabla u^3|^2$ is also good.
By Young's inequality, the  claim (5.10) is proved.

Using (5.9) and (5.10), a standard iteration (cf. \cite {G},
Chapter III, Lemma 2.1) yields that for all $\rho \leq R_0$, one
has
$$\int_{Q_{\rho}} |\nabla u|^2\,dz\leq C\rho^{3+3\b }, \tag 5.11 $$
where $C$ depends on $R_0$. An iteration by (5.9) and (5.10)
yields that for any $\sigma <1$,
$$\int_{Q_{\rho}} |\nabla u|^2\,dz\leq C\rho^{3+2\sigma  }.$$
Using (5.2) and (5.8) yields
$$\split
\int_{Q_{\rho}} |\nabla u -(\nabla u)_{\rho} |^2\,dz&\leq C \left
(\frac {\rho} R\right )^{7}\int_{Q_R} |\nabla
u-(\nabla u)_R|^2\,dz +CR^{\b}\int_{Q_R}|\nabla u|^2\,dz\\
&\leq C \left (\frac {\rho} R\right )^{7}\int_{Q_R} |\nabla
u-(\nabla u)_R|^2\,dz +CR^{3+2\sigma +\b} .\endsplit
$$
Choose $\sigma$ sufficiently close to 1 so that $2\sigma +\b
>2$.
 Then, for all $\rho\leq \frac R2$, we have
$$\int_{Q_{\rho}} |\nabla u -(\nabla u)_{\rho} |^2\,dz\leq
C\rho^{5+2 \sigma_1}
$$
for some $\sigma_1$ with $0<\sigma_1<1$. This implies $\nabla u\in
C^{1,\sigma_1}_{loc}$ and then $\nabla u\in C^{1,\b}$ (cf \cite
{GS}).\qed
\enddemo

\proclaim {Proposition 15} Let $(u, v)$ be  a weak solution of
(1.5)-(1.7) in $\R^2\times [0,T]$ and assume that $u$ is H\"older
continuous in $\R^2\times [0,T)$.
 Let $\tau$ be any positive constant. For
$t\in [\tau, T]$, we have
$$\int_{\R^2} |\nabla^2 u(x,t) |^2+|\nabla v (x,t)|^2\,dx\leq C\,\tau^{-1} (1 +T
R^{-2}).
$$
Then, $(u,v)$ is smooth in $\R^2\times (0,T)$.
\endproclaim
\demo{Proof} By Sobolev's embedding Theorem, we have
$$\int_{B_1(x_0)} |\nabla u(x,t) |^p+| v (x,t)|^p\,dx
$$
for any $p>1$ and for $x_0\in\R^2$ and $t>\tau$. By a similar way
to one in Lemma 5, we can show that $u$ is H\"older continuous in
$\R^2\times [\tau, T]$.

To get the higher order regularity, we rewrite (1.7) as
$$u^i_t - \nabla_{\a} \left [W_{p_{\a}^i}
(u,\nabla u)  \right ] = - u^ku^i \nabla_{\a} \left
[V_{p_{\a}^k}(u,\nabla u) \right ]-(v\cdot \nabla  )u^i +\tilde
B(u,\nabla u),\tag 5.12
$$
where $\tilde B(u, \nabla u)$ is given by
$$\split  \tilde B(u, \nabla u)=&-W_{u^i}(u,\nabla
u)+W_{u^k}(u,\nabla u)u^k u^i+W_{p_{\a}^k}(u,\nabla
u)\nabla_{\a}u^k u^i \\
&+V_{p_{\a}^k}(u,\nabla u) u^k\nabla_{\a}u^i -\nabla_{\a} \left [
u^ku^i \right ]V_{p_{\a}^k}(u,\nabla u).
\endsplit
$$
Since $W(u,p)$ is quadratic and convex in $p$, we can write
$$W_{p^i_{\a}}(u,\nabla u)=a^{ij}_{\a\b} (u)\nabla_{\a}u^j. $$
 Since  $u$ is uniformly H\"older continuous, the
left-hand term of (5.12) is a parabolic operator. Let $\xi (x)$ be
a cut-off function in $B_R(x_0)$ and let $\tau\in C^{\infty}(\R,
\R)$ be a function depends only on $t$ with $0\leq \tau\leq 1$,
$\tau\equiv 1 $ on $[t_0-\frac 1 4 R^2, t_0]$ and $\tau\equiv 0$
on $(-\infty , t_0-R^2)$ and $|\partial_t \tau|\leq C/R^2$. Set
$\phi =\tau \xi$. Multiplying (5.12) by $\phi$, we have
$$\split
&\quad (u\phi)^i_t - \nabla_{\a} \left [ a^{ij}_{\a\b}
(u)\nabla_{\a}
(u^j \phi) \right ]- u^i \phi_t \\
&= - u^ku^i \nabla_{\a} \left [ V_{p_{\a}^k}(u,\nabla u) \right
]\phi -[(v\cdot \nabla  )u^i +\tilde B(u,\nabla u)] \phi.
\endsplit\tag 5.13
$$
By the assumption, we have
$$(v\cdot \nabla )u\in L^p (Q_{R}(x_0)) ,\quad |\nabla u|^2\in  L^p (Q_{R}(x_0) )\quad \forall p>1 $$
But the first term on the righthand of (5.13) is not a `good'
term, which need more analysis. Using the fact that $|u|=1$, we
have
$$u^3\nabla^2_{\a\b} u^3=-(\nabla_{\b}u\cdot\nabla u +
u^1\nabla_{\a\b}^2 u^1+ u^2\nabla_{\a\b}^2 u^2),\quad u^3 u_t^3=-(
u^1 u_t^1+ u^2  u_t^2).
$$
Without loss of generality, we regard the solution in $\R^3$. By a
rotation, we assume
$$u(x_0,t_0)=(0,0,1).
$$
Since $u$ is H\"older continuous, there exists a small $R$ such
that
$$|u(x,t)-u(x_0,t_0)|\leq \varepsilon
$$
for a sufficiently small constant $\varepsilon >0$. Therefore
$$|\nabla ^2 u^3|\leq  C|\nabla u|^2+ 2\varepsilon ( |\nabla ^2
u^1|+|\nabla ^2 u^2|)
$$

Apply the classical $L^p$-estimate of parabolic systems  (c.f.
\cite {Ei}, \cite {LSU}) to (5.13) for $i=1,2$, we have
$$\split
\|\tilde u_t\phi\|_{L^p(Q_R(x_0))} +\|\nabla^2(\tilde u\phi
)\|_{L^p(Q_R(x_0))} &\leq C\|\phi \nabla^2  u^3  \|_{L^p}+
C\varepsilon \|\phi
\nabla^2u\|_{L^p(Q_R(x_0))}\\
&+C(\|u\|_{L^{2p}(Q_R(x_0))}+\|v\|_{L^{2p}(Q_R(x_0))} +1),
\endsplit
$$
where $\tilde u =(u^1,u^2)$. Choosing $\varepsilon$ sufficiently
small, we obtain
$$\|u_t\phi\|_{L^p(Q_R(x_0))} +\|\nabla^2(u\phi )\|_{L^p(Q_R(x_0))}\leq C.
$$

To estimate $v$ in (1.5), it follows from H\"older's inequality
that
$$\int_{\R^2\times [\tau ,T]}|(v\cdot \nabla )v|^p\,dx\leq
\left (\int_{\R^2\times [\tau ,T]}|\nabla v|^4\,dx\,dt\right
)^{p/4} \left (\int_{\R^2\times [\tau ,T]}|v|^{\frac
4{4-p}}\,dx\,dt\right )^{\frac {4-p} 4}
$$
for any $p$ with $3<p<4$.
 By the $L^p$-estimate of Stoke's
operator (e.g. \cite {So}),  $v_t$ and $\nabla^2 v$ are in $L^p$
for $3<p<4$. This implies that $v$ is H\"older continuous.

Differentiating in $x_l$ in (5.12), we have
$$\split
&(\nabla_{x_l}u^i)_t - \nabla_{\a} \left [ a^{ij}_{\a\b}
(u)\nabla_{\a} (\nabla_{x_l}u^j)  \right ]\\
& = - u^ku^i \nabla_{\a} \left [V_{p_{\a}^k}(u,\nabla \nabla_{x_l}
u) \right ]+v\#\nabla^2u+ \nabla v\#\nabla u+ \nabla u\#\nabla^2 u
.
\endsplit
$$
By  applying the $L^p$-theory, a similar argument yields that
$\nabla u$ is uniformly continuous.
 Then, a
standard bootstrap method implies that $(u,v)$ are smooth.\qed
\enddemo
\bigskip\noindent {\it Acknowledgements:} The research  of the
first author was supported by the Australian Research Council
grant DP0985624. The research of the second author is supported
partially by Zheng Ge Ru Foundation, Hong Kong RGC Earmarked
Research Grant CUHK4042/08P, and a Focus Area Grant from The
Chinese University of Hong Kong. A part of the work was done when
Hong visited the Chinese University of Hong Kong in July of 2008
and in December of 2009.


\Refs\widestnumber\key {[CKMS]}

\ref \key {AL} \by F. J. Almgren and E. H. Lieb\paper Sigularities
of energy minimizing maps from the ball to the sphere: Examples,
counterexamples, and bounds\jour Ann. Math. \vol 128 \yr 1988
\pages 483--530
\endref


\ref \key {Am} \by H. Amann\paper Quasilinear  Parabolic systems
under nonlinear boundary conditions\jour Arch. Rational Mech.
Anal. \vol 92\yr 1986\pages 153-–192
\endref

\ref \key {BCLP} \by P. Bauman, M. Calderer, C. Liu and D.
Phillips\paper The Phase Transition between Chiral Nematic and
Smectic A. Liquid Crystals \jour Arch. Rational Mech. Anal. \vol
165\yr  2002\pages  161-–186
\endref

\ref \key {BBC}\by F. Bethuel, H. Brezis and J. M. Coron \paper
Relaxed energies for harmonic maps\jour In variational methods,
edited by Berestycki, Coron, Ekeland, Birkh\"auser, Basel \yr
1990\pages 37--52\endref


\ref \key {CKN}\by L. Caffarelli, R. Kohn and L. Nirenberg \paper
Partial regularity of suitable weak solutions of Navier-Stokes
euqations\jour Comm. Pure Appl. Math.\vol 35 \yr 1982\pages
771--831
\endref

\ref \key {CS} \by Y. Chen and M. Struwe \paper  Existence and
partial regular results for the heat flow for harmonic maps\jour
Math. Z. \vol  201 \yr 1989 \pages 83--103\endref

\ref \key {G}\by M. Giaquinta\book Multiple integrals in the
calculus of variations and nonlinear elliptic systems\publ
Princeton Univ. Press \yr 1983 \endref

\ref \key {GMS1}\by M. Giaquinta, G. Modica and J. Soucek\paper
The Dirichlet energy of mappings with values into the sphere \jour
Manuscripta Math.\vol 65\yr 1989\pages 489--507\endref


\ref \key {GMS2}\by M. Giaquinta, G. Modica and J. Soucek\paper
Liquid crystals: relaxed energies, dipoles, singular lines and
singular Points \jour Ann. Scuola Norm. Sup. Pisa (3)\vol 17\yr
1990\pages 415--437\endref


\ref \key {GMS3}\by M. Giaquinta, G. Modica and J. Soucek\book
Cartesian currents in the calculus of variations, part II,
Variational integrals \publ  A series of modern surveys in
mathematics, 38, Springer-Verlag\yr 1998 \endref

\ref \key {GS}\by M. Giaquinta and M. Struwe \paper On the partial
regularity  weak solutions of non-linear parabolic systems\jour
Math. Z\vol  179\yr 1982\pages 437--451\endref

\ref \key {ES}\by J. Eells and J. H. Sampson
 \paper Harmonic mappings of Riemannian manifolds\jour Amer. J. Math.\vol  86 \yr 1964\pages
109--160 \endref


\ref \key {Ei}\by S. Eidel'man\book Parabolic systems \publ North
Holland Publishing\yr 1969
\endref

\ref \key {Er} \by J. Ericksen  \book  Equilibrium Theory of
Liquid Crystals \publ Academic Press, New York  \yr 1976\endref

\ref \key {HKL1} \by R. Hardt, D. Kinderlehrer and F.-H. Lin\paper
Existence and partial regularity of static Liquid Crystal
Configurations \jour Comm.  Math. Phys.\vol 105 \yr 1986\pages
547--570\endref

\ref \key {HKL2} \by R. Hardt, D. Kinderlehrer and F.-H. Lin\paper
Stable defects of minimizers of constrained variational principles
\jour Ann. Inst. Henri Poincar\'e, Analyse non lin\'eaire\vol 5
\yr 1988\pages 297--322\endref

\ref \key {Ho1}\by M.-C. Hong \paper Partial regularity of weak
solutions of the Liquid Crystal equilibrium system \jour Indiana
Univ. Math. J. \vol  53  \yr 2004 \pages1401-1414\endref

\ref \key {Ho2}\by M.-C. Hong \paper  Existence of infinitely many
equilibrium configurations of the Liquid Crystal system
prescribing the same non-constant boundary value\jour Pacific
Journal of Mathematics \vol  232 \yr 2007 \pages 177--206\endref

\ref \key {Ho3}\by M.-C. Hong \paper  Global existence of
solutions of the simplified Ericksen-Leslie system  in dimension
two \jour To appear in Calc. Var. $\&$ PDEs (DOI
10.1007/s00526-010-0331-5)\endref


\ref \key  {Hu}\by N. Hungerb\"{u}hler, $m$-harmonic flow\jour
Ann. Scuola Norm. Sup. Pisa Cl. Sci. \vol (4) XXIV \yr 1997 \pages
593--631\endref


\ref \key {K}\by M. Kleman
 \book  Points, Lines and Walls
\publ John Wiley \& Son, New Year\yr 1983\endref


\ref \key {LSU} \by O. A. Ladyzhenskaya, V. A. Solonnikov and N.
N. Ural'ceva\book Linear and qusilinear equations of parabolic
type\publ  Tanslations of Mathematical Monographs 23. Providence,
Rhode Island: American Mathematical Society\yr 1968\endref

\ref \key {Le}\by F. Leslie\book Theory of flow phenomenon in
liquid crystal   \publ (vol 4) Brown (Ed.) A. P. New York \yr
1979\pages 1--81
\endref


\ref \key {L1}\by F.-H. Lin \paper Nonlinear theory of defects in
nematic liquid crystals: Phase transition and flow phenomena\jour
Comm. Pure Appl. Math.\vol 42 \yr 1989\pages 789--814
\endref

\ref \key {L2}\by F.-H. Lin \paper A new proof of the
Caffarelli-Kohn-Nirenberg theorem \jour Comm. Pure Appl. Math.\vol
  51  \yr 1998\pages  241--257
\endref

\ref \key {LLW}\by  F.-H. Lin, J. Lin and C. Wang \paper Liquid
crystal flow in two
 dimension \jour
 To appear in  Arch. Rational Mech.
Anal.\endref

 \ref \key {LL1}\by F.-H. Lin and C. Liu\paper Nonparabolic
dissipative systems modelling the flow of liquid cystals\jour
Comm. Pure Appl. Math.\vol 48 \yr 1995\pages 501--537
\endref

\ref \key {LL2}\by F.-H. Lin and C. Liu\paper Existence of
solutions for the Ericksen-Leslie System\jour Arch. Rational Mech.
Anal.\vol 154 \yr 2000\pages 135--156
\endref

\ref \key {LP}\by F.-H. Lin and X.-B. Pan\paper Magnetic
field-induced instabilities in liquid crystals \jour SIAM J. Math.
Anal \vol 38 \yr 2007\pages 1588--1612
\endref

\ref \key {So} \by  V. A. Solonnikov   \paper $L_p$-estimates for
solutions to the initial boundary-value problem for the
generalized Stokes system in a bounded domain.  \jour  J. Math.
Sci.\vol  105\yr  2001\pages  2448--2484 \endref

 \ref \key {St1} \by M. Struwe  \paper
On the evolution of harmonic maps of Riemannian surfaces\jour
Commun. Math. Helv.\vol 60\yr 1985\pages 558--581 \endref

\ref \key {St2} \by M. Struwe  \paper The existence of surfaces of
constant mean curvature with free boundaries \jour Acta Math.\vol
160\yr 1988\pages 19-64 \endref


\ref \key {Sc} \by V. Scheffer \paper Hausdorff measure and the
Navier-Sokes equations\jour Comm. Math. Phys.\vol 61\yr 1977\pages
97--112
\endref

\ref \key {TX} \by G. Tian and Z. Xin  \paper Gradient estimation
on Navier-Stokes equations\jour Comm. Anal. Geom. \vol  7\yr 1999
\pages 221--257
\endref


\endRefs
\enddocument
\end